\documentclass{amsart}
\usepackage{amssymb,euscript,amsmath, mathrsfs}
\usepackage[dvips]{graphicx}
\usepackage[dvips]{color}
\usepackage{epsfig}

\newcounter{ENUM}
\newcommand{\itm}{\item}
\newenvironment{ilist}{\renewcommand{\theENUM}{\roman{ENUM}}\renewcommand{\itm}{\addtocounter{ENUM}{1}\item[(\theENUM)]}\begin{itemize}\setcounter{ENUM}{0}}{\end{itemize}}
\newenvironment{alist}[1][0]{\renewcommand{\theENUM}{\alph{ENUM}}\renewcommand{\itm}{\addtocounter{ENUM}{1}\item[\theENUM)]}\begin{itemize}\setcounter{ENUM}{#1}}{\end{itemize}}

\newcommand{\nHat}[1]{\widehat{#1}}                        

\def\S{{\mathfrak S}}
\def\e{{\mathfrak e}}
\def\bf{{\mathbf f}}
\def\bg{{\mathbf g}}
\def\be{{\mathbf e}}

\def\x{{\mathbf x}}
\def\y{{\mathbf y}}
\def\z{{\mathbf z}}

\def\v{{\mathbf v}}
\def\u{{\mathbf u}}
\def\w{{\mathbf w}}
\def\s{{\mathbf s}}
\def\1{{\mathbf 1}}

\def\Z{{\mathbb Z}}
\def\N{{\mathbb N}}

\def\Q{{\mathbb Q}}
\def\R{{\mathbb R}}

\def\H{{\mathbb H}}

\def\F{{\mathscr F}}

\def\cP{{\mathcal P}}

\def\f{{\mathfrak f}}

\def\sgn{\operatorname{sign}}

\def\conv{\mathrm{conv}}
\def\aff{\mathrm{aff}}
\def\lin{\mathrm{lin}}

\def\vol{\mathrm{Vol}}
\def\svol{\mathrm{SVol}}
\def\vert{\mathrm{vert}}
\def\Int{\mathrm{Int}}
\def\tcone{\mathrm{tcone}}
\def\fcone{\mathrm{fcone}}

\newtheorem{thm}{Theorem}[section]
\newtheorem{prop}[thm]{Proposition}
\newtheorem{lem}[thm]{Lemma}
\newtheorem{cor}[thm]{Corollary}

\theoremstyle{definition}
\newtheorem{defn}[thm]{Definition}

\newtheorem{ex}[thm]{Example}

\theoremstyle{remark}

\newtheorem{rem}[thm]{Remark}

\numberwithin{equation}{section}

\address{Department of Mathematics, University of California, One Shields Avenue, Davis, California 95616.}
\email{fuliu@math.ucdavis.edu}

\begin{document}
\title{Higher integrality conditions, volumes and Ehrhart polynomials}
\author{Fu Liu}\thanks{The author is partially supported by the National Security Agency under Grant No. H98230-09-1-0029.}
\begin{abstract}
A polytope is integral if all of its vertices are lattice points. The constant term of the Ehrhart polynomial of an integral polytope is known to be $1$. In previous work, we showed that the coefficients of the Ehrhart polynomial of a lattice-face polytope are volumes of projections of the polytope. We generalize both results by introducing a notion of $k$-integral polytopes, where $0$-integral is equivalent to integral. We show that the Ehrhart polynomial of a $k$-integral polytope $P$ has the properties that the coefficients in degrees less than or equal to $k$ are determined by a projection of $P$, and the coefficients in higher degrees are determined by slices of $P$. A key step of the proof is that under certain generality conditions, the volume of a polytope is equal to the sum of volumes of slices of the polytope.
\end{abstract}

\maketitle

\section{Introduction}

Let $V$ be a finite-dimensional real vector space. A subset $\Lambda \subset V$ is called a {\it lattice of $V$} if $\Lambda$ is a discrete additive subgroup of $V$ and spans $V.$ Let $\bar{\Lambda}$ be the rational vector space generated by $\Lambda.$ Any point in $\Lambda$ is called a {\it lattice point}, and any point in $\bar{\Lambda}$ is called a {\it rational point}. 

A subset $\Gamma \subset \Lambda$ is a {\it sublattice of $\Lambda$} if $\Gamma$ is a subgroup of $\Lambda.$ 
We say a sublattice $\Gamma$ is {\it of rank $r$} if $\aff(\Gamma)$ has dimension $r.$

We assume familiarity with basic definitions of polytopes as presented in \cite{zie}. An {\it integral} polytope is a convex polytope whose vertices are all lattice points. A {\it rational} polytope is a convex polytope whose vertices are all rational points. 
For any sublattice $\Gamma \subset \Lambda,$ let $U$ be the subspace of $V$ spanned by $\Gamma.$ For any polytope $P$ lying in an affine space that is a translation of $U,$ we define the {\it volume of $P$ normalized to the lattice $\Gamma$} to be the integral
\begin{equation}\label{defVol}
\vol_{\Gamma}(P) := \int_P 1 \ {d \ \Gamma},
\end{equation}
where $d \ \Gamma$ is the canonical Lebesgue measure defined by the lattice $\Gamma.$ 

For any point $\y \in V$ and any set $S \subset V,$ we denote by $S + \y$ the set $\{ \s + \y \ | \ \s \in S\}$ obtained from translating $S$ by $\y.$ For any affine space $W \subset V,$ we denote by $\lin(W)$ the translation of $W$ to the origin. It is easy to see that $W = \lin(W) + \y,$ for some (equivalently all) $\y$ in $W.$

If $P$ is a rational polytope in $V,$  we denote by $\lin(P)$ the translation of $\aff(P)$ to the origin and denote by $\Lambda_{\lin(P)}$ the lattice $\lin(P) \cap \Lambda.$ One checks $\Lambda_{\lin(P)}$ spans $\lin(P).$ Thus, we can define $\vol_{\Lambda_{\lin(P)}}(P),$ which we often refer as the {\it normalized volume} of $P.$ 

Throughout the paper, we will assume $V$ has dimension $D,$ (so $\Lambda$ is of rank $D,$) and fix a basis $\e = (\be^1, \be^2, \dots, \be^D)$ of the lattice $\Lambda.$ Then we can represent each point in $V$ by its coordinates with respect to $\e,$ namely, for any point $\sum_{i=1}^D c_i \be^i \in V,$ we say its {\it coordinates (with respect to $\e$)} are $(c_1, \dots, c_D).$ This gives us an isomorphism between $V$ and $\R^D$, and $\Lambda$ and $\Z^D.$  For convenience, we denote by $\Lambda_k := \langle \be^1, \dots, \be^k \rangle$ the lattice generated by $\be^1, \dots, \be^k$, $\Lambda^k := \langle \be^{k+1}, \dots,  \be^D \rangle$ the lattice generated by  $\be^{k+1}, \dots, \be^D,$ $V_k$ the $k$-dimensional subspace of $V$ spanned by $\be^1, \dots, \be^k,$ and $V^k$ the $(d-k)$-dimensional subspace of $V$ spanned by $\be^{k+1}, \dots, \be^{D}.$ Let $\pi^{(k)}: V \to V_{D-k}$ be the projection that maps $\sum_{i=1}^D c_i \be^i$ to $\sum_{i=1}^{D-k} c_i \be^i,$ or equivalently the map that forgets the last $k$ coordinates.

For any polytope $P$ in $V,$ we denote by
$$i(P) := \#(P \cap \Lambda)$$
the number of lattice points inside $P.$ Furthermore, for positive integer $m \in \N,$ we denote by $i(P,m)$ the number of lattice points in $m P,$ where $m P = \{ m x \ | \ x \in P \}$
is the {\it $m$th dilation} of $P.$ Eug\`{e}ne Ehrhart \cite{Ehrhart} discovered that for any $d$-dimensional integral polytope, $i(P,m)$ is a polynomial of degree $d$ in $m.$ Thus, we
call $i(P,m)$ the {\it Ehrhart polynomial} of $P.$ Although Ehrhart's
theory was developed in the 1960's, we still do not have a very good understanding of
the coefficients of Ehrhart polynomials for general polytopes except
that the leading, second and last coefficients of $i(P,m)$ are the
normalized volume of $P$, one half of the normalized volume of the
boundary of $P$ and $1,$ respectively. In \cite{cyclic}, the author showed that for any $d$-dimensional integral cyclic polytope $P,$ we have that \begin{equation}\label{ques}
i(P, m) =  \vol_\Lambda(mP) + i(\pi(P), m) = \sum_{j=0}^D
\vol_{\Lambda_j}(\pi^{(D-j)}(P)) m^j,
\end{equation}
where $D = d.$
In \cite{lattice-face, note-lattice-face}, the author generalized the family of integral cyclic polytope to a bigger family of integral polytopes, {\it
lattice-face} polytopes, and showed that their Ehrhart polynomials also satisfy (\ref{ques}).

The motivation of this paper is to prove a conjecture given in \cite[Conjecture 8.5]{lattice-face}, which can be viewed as generalizations of both the theorem for lattice-face polytopes shown in \cite{lattice-face, note-lattice-face} and the fact that the constant term of the Ehrhart polynomial of an integral polytope is $1.$ In fact, we will prove a stronger version of this conjecture. We need the following definitions before we can state the result.

\begin{defn}
Any affine space $U$ in $V$ is {\it integral with respect to $\e$} if
$$\pi^{(D-\dim(U))}(U \cap \Lambda) = \Lambda_{\dim(U)}.$$
A polytope is {\it affinely integral with respect to $\e$} if its affine hull is integral with respect to $\e$.

Suppose $0 \le k \le d \le D.$ A $d$-dimensional polytope in $V$ is {\it $k$-integral with respect to $\e$} if any face of $P$ of dimension less than or equal to $k$ is affinely integral with respect to $\e.$

In particular, when $k = d,$ we call $P$ a {\it fully integral} polytope with respect to $\e.$
\end{defn}

Note that the definition of $0$-integral polytopes is the same as that of integral polytopes. Thus, any $k$-integral polytope is an integral polytope. Since we fix the lattice basis $\e$ in the paper, we often omit ``with respect to $\e$'' unless it is not clear from the context.

\begin{defn}
Let $S$ be a set in $V$ and $0 \le k \le d.$ For any $\y \in \pi^{(D-k)}(V),$ we define $\pi_{D-k}(\y, S)$ to be the intersection of $S$ with the inverse image of $\y$ under $\pi^{(D-k)},$ and we call it {\it the slice of $S$ over $\y$.}

\end{defn}

Our first main theorem is the following.

\begin{thm}\label{mainEhr}
Suppose $0 \le k \le d \le D$ and $P  \subset V$ is a $k$-integral $d$-dimensional polytope. Then the Ehrhart polynomial of $P$ is given by
\begin{eqnarray*}
i(P,m) &\!=\!\!& m^k \left( \sum_{\y \in \pi^{(D-k)}(P) \cap \Lambda} i(\pi_{D-k}(\y,P), m)-1 \right)  + i(\pi^{(D-k)}(P), m) \\
&\! =\!\!& 
m^k  \left( \sum_{\y \in \pi^{(D-k)}(P) \cap \Lambda} i(\pi_{D-k}(\y,P), m)-1 \right) + \sum_{j=0}^k \vol_{\Lambda_j} (\pi^{(D-j)}(P)) m^j.
\end{eqnarray*}
\end{thm}

We will see in Corollary \ref{sliceInt} that each $\pi_{D-k}(\y, P)$ appearing in the above formula is an integral polytope, so the constant term of $i(\pi_{D-k}(\y,P), m)$ is $1.$ Therefore, Theorem \ref{mainEhr} says that the coefficient of $m^j$ in the Ehrhart polynomial of a $k$-integral $d$-dimensional polytope $P$ is 
$\vol_{\Lambda_j} (\pi^{(D-j)}(P))$ for $0 \le j \le k$, and is the coefficient of $m^{j-k}$ in the sum of the Ehrhart polynomials of the slices of $P$ over lattice points in $\pi^{(D-k)}(P)$ for $k+1 \le j \le d.$

One consequence of Theorem \ref{mainEhr} is that the coefficients of $i(P,m)$ are all positive for any $d$-dimensional polytope that is $(d-2)$-integral (with respect to some lattice basis). This provides one possible way to prove positivity conjectures on coefficients of Ehrhart polynomials of special families of polytopes. 
It was recently conjectured that the coefficients of the Ehrhart polynomials of both Birkhoff polytopes and matroid polytopes are positive. One could approach these conjectures by proving these polytopes are $(d-2)$-integral, or by showing there exist subdivisions of these polytopes into $(d-2)$-integral polytopes (not necessarily with respect to a fixed lattice basis) satisfying certain conditions.

Recall that the leading coefficient of the Ehrhart polynomial of a polytope is the normalized volume of the polytope. If we only look at the leading coefficients in Theorem \ref{mainEhr}, we get that the normalized volume of $P$ is equal to the sum of the normalized volume of slices of $P$ over lattice points in $\pi^{(D-k)}(P),$ provided that $P$ is $k$-integral. In fact, we can relax the condition of $k$-integral by introducing the following definition.

\begin{defn}
Any affine space $U$ in $V$ is {\it in general position with respect to $\e$} if
$$\pi^{(D-\dim(U))}(U) = V_{\dim(U)}.$$
A polytope is {\it in affinely general position with respect to $\e$} if its affine hull is in general position with respect to $\e$.

Suppose $0 \le k \le d \le D.$ A $d$-dimensional polytope in $V$ is {\it in $k$-general position with respect to $\e$} if any face of $P$ of dimension less than or equal to $k$ is in affinely general position with respect to $\e.$

In particular, when $k = d,$ we say $P$ is {\it in fully general position  with respect to $\e$}.
\end{defn}

Again, we often omit ``with respect to $\e$'' if there is no possibility of confusion.


Our result on volumes and slices is the following.

\begin{thm}\label{mainVol}
Suppose $0 < k < d \le D$ and $P \subset V$ is a $d$-dimensional $(k-1)$-integral polytope in $k$-general position. Let $\Lambda_P := \aff(P) \cap \Lambda.$ Then the normalized volume of $P$ is given by
\begin{equation}\label{VolFor0}
\vol_{\Lambda_{\lin(P)}}(P) = \sum_{\y \in \pi^{(D-k)}(\Lambda_P)} \vol_{\Lambda^k \cap \Lambda_{\lin(P)}}(\pi_{D-k}(\y,P)).
\end{equation}
\end{thm}

We remark that Theorem \ref{mainVol} is trivially true when $k =0$ because the right hand side of \eqref{VolFor0} becomes $\vol_{\lin(P)}(P)$. However, the theorem does not hold for $k=d.$ The most interesting case of Theorem \ref{mainVol} is when $k=1.$ If $k=1,$ the condition becomes that $P$ is an integral polytope and is in $1$-general position. However, for any integral polytope, it is always in $1$-general position with respect to some basis. Hence, formula \eqref{VolFor0} with $k=1$ can be applied to any integral polytope. Moreover, any rational polytope can be dilated to an integral polytope, so one can obtain a formula modified from \eqref{VolFor0} with $k=1$ to calculate volumes of rational polytopes.

Although Theorem \ref{mainVol} can be viewed as a corollary to Theorem \ref{mainEhr}, it actually serves as a key step in proving Theorem \ref{mainEhr}. In fact, most of the paper will be dedicated to the proof of Theorem \ref{mainVol}.

The plan of this paper is as follows: In Section \ref{exthms}, we show examples of our main theorems. In Section \ref{prelim}, we introduce basic definitions and lemmas and reduce Theorem \ref{mainVol} to Proposition \ref{mainVol1}. In Sections \ref{preservation}, \ref{reduction} and \ref{special}, we prove a special case of Proposition \ref{mainVol1} when $k=1$ and $P$ is a simplex.  Our method is to reduce the problem to the full dimensional, fully general position case, and then to carry out the calculation directly using techniques developed in \cite{lattice-face} involving signed decomposition, determinants, and power sums.
In Section \ref{propertySP}, we discuss the properties of slices and projections of a polytope. Using these, we complete the proof of Proposition \ref{mainVol1} and thus the proof of Theorem \ref{mainVol} in Section \ref{pfVol}. Finally, in Section \ref{pfEhr}, we prove Theorem \ref{mainEhr} by using Theorem \ref{mainVol} and the local formula  of Berline and Vergne relating the number of lattice points to volumes of faces for rational polytopes \cite{berline-vergne}.

\section{Examples of the theorems}\label{exthms}
In this section, we give examples of the main theorems and also show that the conditions in Theorem \ref{mainVol} are in fact necessary. Throughout the section, we assume that $D=d,$ $\Lambda = \Z^d$ with the standard basis, and $V = \R^d.$ 

\begin{ex}[Example of Theorem \ref{mainEhr}]\label{exEhr}
Consider the $3$-dimensional polytope $$P = \conv \{(0,0,0), (4,0,0), (3,6,0), (2,2,2)\} \subset \R^3.$$ One checks that $P$ is $1$-integral.

\begin{figure}
\centerline{
\mbox{\includegraphics[scale=0.4, viewport=90 125 420 350,clip]{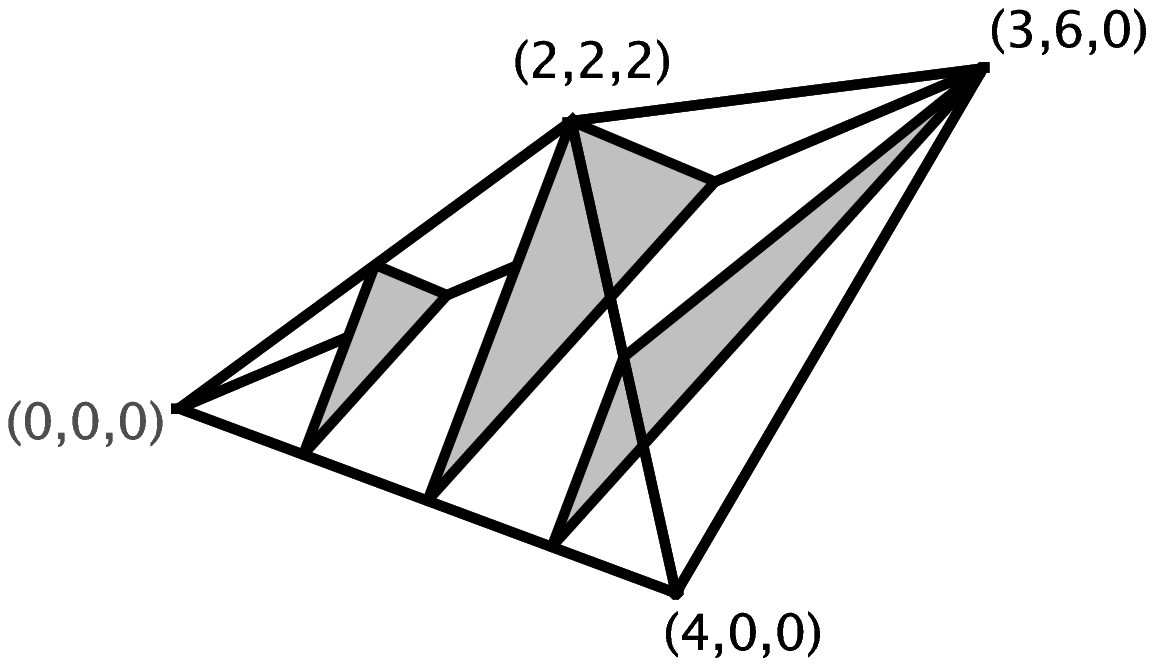}} \ \ \ \ \
\mbox{\includegraphics[scale=0.4, viewport=90 125 420 350,clip]{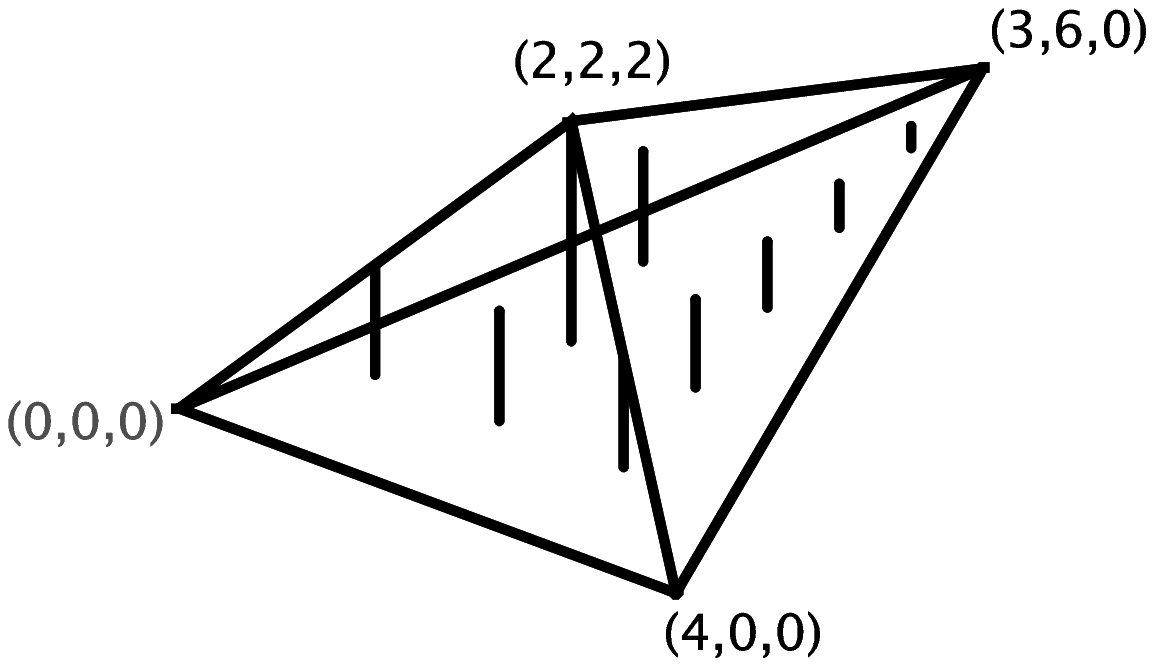}}
}
\caption{Slices of a polytope}\label{slicepic}
\end{figure}

$\pi^{(2)}(P) = [0,4].$ The lattice points in $\pi^{(2)}(P)$ are $0,1,2,3$ and $4.$ In the picture on the left side of Figure \ref{slicepic}, the three shaded triangles are the slices of $P$ over the lattice points $1, 2$ and $3.$ The other two slices of $P$ over lattice points are the single points $\pi_{2}(0, P) = (0,0,0)$ and $\pi_{2}(4,P) = (4,0,0).$ We calculate the Ehrhart polynomial of each slice and get $i(\pi_2(0,P),m) = 1, i(\pi_2(1,P),m) = m^2 + 2m +1, i(\pi_2(2,P),m)=4m^2+4m+1, i(\pi_2(3,P),m) = 3m^2 + 4m + 1$ and  $i(\pi_2(4,P),m) = 1.$ Their sum is 
$${\textbf 8}m^2 + \textbf{10}m + 5.$$
We also have that
$$\vol_{\Lambda_1}(\pi^{(2)}(P)) = {\textbf 4}, \mbox{ and } \vol_{\Lambda_0}(\pi^{(3)}(P))={\textbf 1}$$

Therefore, applying Theorem \ref{mainEhr} with $k=1,$ we conclude that $8, 10, 4$ and $1$ are the coefficients for $m^3, m^2, m^1$ and $m^0$ in $i(P,m)$. Hence, the Ehrhart polynomial of $P$ is given by 
$$i(P,m) = 8m^3 + 10m^2 + 4m +1.$$
\end{ex}

The rest of the examples in this section will be related to Theorem \ref{mainVol}. Since all the polytopes we will consider are full-dimensional, we have that  $\Lambda_{\lin(P)} = \Lambda_P = \Lambda,$ $\pi^{(D-k)}(\Lambda_P) = \Lambda_k$ and $\Lambda^k \cap \Lambda_{\lin(P)} = \Lambda^k.$ Hence, \eqref{VolFor0} can be simplified to
\begin{equation}\label{VolForFullDim}
\vol_{\Lambda}(P) =\sum_{\y \in \Lambda_{d-k}} \vol_{\Lambda^k}(\pi_{d-k}(\y,P)).
\end{equation}

\begin{ex}[Example of Theorem \ref{mainVol}]
Let $P$ be the same polytope as in Example \ref{exEhr}.
Although $P$ is not $2$-integral, it is $1$-integral and in $2$-general position. We will apply Theorem \ref{mainVol} on $P$ for $k=2.$ 
The volume of $P$ (normalized to $\Lambda$) is $8.$ We will show the right hand side of \eqref{VolFor0}, or equivalently the right hand side of \eqref{VolForFullDim}, gives $8$ as well.

$\pi^{(1)}(P) = \conv \{(0,0),(4,0),(3,6) \}.$ There are $9$ lattice points in the interior of $\pi^{(1)}(P)$ and $8$ lattice points on the boundary. In the picture on the right side of Figure \ref{slicepic}, the $9$ line segments are the slices of $P$ over the lattice points in the interior of $\pi^{(1)}(P)$. For example, the left-most line segment is the slice of $P$ over $(1,1):$
$$\pi_1((1,1), P) = \conv\{(1,1,0), (1,1,1)\}.$$
The slices over the boundary lattice points are all single points, thus have $0$ $\Lambda^2$-volume. Hence, we only calculate the $\Lambda^2$-volume of the slices over interior lattice points: 
$$\begin{array}{lll}
\vol_{\Lambda^2}(\pi_1((1,1), P)) =  1,  & \vol_{\Lambda^2}(\pi_1((2,1), P))=1,  &\vol_{\Lambda^2}(\pi_1((2,2), P))=2, \\
\vol_{\Lambda^2}(\pi_1((2,3), P))=1, & \vol_{\Lambda^2}(\pi_1((3,1), P))=1, & \vol_{\Lambda^2}(\pi_1((3,2), P))=\frac{4}{5}, \\
\vol_{\Lambda^2}(\pi_1((3,3), P))=\frac{3}{5}, & \vol_{\Lambda^2}(\pi_1((3,4), P))=\frac{2}{5}, & \vol_{\Lambda^2}(\pi_1((3,5), P))=\frac{1}{5}.
 \end{array}$$
The sum is $8,$ as desired.
\end{ex}

The hypothesis in Theorem \ref{mainVol} is that $P$ is $(k-1)$-integral and in $k$-general position. One might wonder whether it is possible to relax the condition and still get the same result \eqref{VolFor0}. There are two natural ways to relax the condition: 1) only requiring $P$ to be $(k-1)$-integral; 2) for $k \ge2,$ requiring $P$ to be $(k-2)$-integral and in $k$-general position. In the next two examples, we show that  \eqref{VolFor0} or equivalently \eqref{VolForFullDim} does not hold if we relax the hypothesis in Theorem \ref{mainVol} to either 1) or 2).

\begin{ex}
Consider $P$ to be the unit square in $\R^2.$ Then $P$ is integral, but not in $1$-general position, because it has edges $e_1 = \conv\{(0,0),(0,1)\}$ and $e_2 = \conv\{(1,0),(1,1)\}$ that are not in affinely general position. 

We have that $\vol_{\Lambda}(P) = 1$ but $\sum_{\y \in \Lambda_1} \vol_{\Lambda^1}(\pi_{1}(\y,P)) = \vol_{\Lambda^1}(e_1) + \vol_{\Lambda^1}(e_2) =2.$

Therefore, Theorem \ref{mainVol} does not hold for $(k-1)$-integral polytopes that are not in $k$-general position.
\end{ex}

\begin{ex}
Consider the polytope $$P = \conv \{(0,0,0), (4,0,0), (3,3,0), (2,1,5)\}.$$ We check that $P$ is an integral polytope that is in $2$-general position. Since the edges $\conv\{(0,0,0),(2,1,5)\}$ and $\conv\{(4,0,0),(2,15)\}$ are not affinely integral, $P$ is not $1$-integral. 

We calculate $\vol_{\Lambda}(P) = 10$, but $\sum_{\y \in \Lambda_2} \vol_{\Lambda^2}(\pi_{1}(\y,P)) = 8.$ 

We thus conclude that Theorem \ref{mainVol} does not hold for $(k-2)$-integral polytopes that are in $k$-general position.
\end{ex}

\section{Preliminaries}\label{prelim}
In this section, we investigate the basic properties of the conditions and formulas involved in our main theorems. The lemmas and definitions introduced here will be used as the foundation of this paper. 

First, we see that any $k$-integral polytope $P$ is integral. Hence, $\aff(P)$ contains a lattice point, say $\beta.$ If we let $Q := P - \beta,$ it is easy to check that to prove Theorem \ref{mainEhr} and Theorem \ref{mainVol} for $P$ is equivalent to proving them for the polytope $Q.$ Note that $Q$ has the property that $\aff(Q)$ contains the origin, thus is a subspace of $V.$ Therefore, we give the following definition:

\begin{defn}
For any polytope $P \subset V,$ we say $P$ is {\it central} if $\aff(P)$ contains the origin, or equivalently, $\aff(P) = \lin(P)$ is a subspace of $V.$
\end{defn}

Based on above discussion, we conclude:
\begin{lem}\label{cent}
It is sufficient to prove Theorem \ref{mainEhr} and Theorem \ref{mainVol} with the assumption that $P$ is central.
\end{lem}

Lemma \ref{cent} is helpful for proving Theorem \ref{mainVol} in particular. If $P$ is central, then $\Lambda_P = \Lambda_{\lin(P)}.$ So formula \eqref{VolFor0} can be simplified to
\begin{equation}\label{VolFor}
\vol_{\Lambda_{P}}(P) = \sum_{\y \in \pi^{(D-k)}(\Lambda_P)} \vol_{\Lambda^k \cap \Lambda_{P}}(\pi_{D-k}(\y,P)).
\end{equation}

\subsection{On lattices}
We know that $\Lambda_k$ is a lattice of $V_k$ and $\Lambda^k$ is a lattice of $V^k.$ In this subsection, we develop a better understanding of $\Lambda_P$ and $\Lambda^k \cap \Lambda_P$.

We say a subspace $U \subset V$ is {\it rational} if it is spanned by a set of vectors with rational coordinates. It is easy to check that a subspace is rational if and only if it is defined by a set of linear equations with rational coefficients.

We state the following three fundamental lemmas about rational subspaces and lattices without proofs. 
\begin{lem}\label{latfund1}
A subspace $U$ of $V$ is rational if and only if there exists a sublattice $\Gamma$ of $\Lambda$ such that $\Gamma$ is a lattice of $U$, if and only if $\Lambda_U := U \cap \Lambda$ is a lattice of $U.$
\end{lem}
\begin{lem}\label{latfund2}
Suppose $\Gamma_i$ is a sublattice of $\Lambda$, and $U_i$ is the subspace of $V$ spanned by $\Gamma_i,$ for $i=1,2.$ Then $\Gamma_1 \cap \Gamma_2$ is a lattice of $U_1 \cap U_2.$
\end{lem}

\begin{lem}\label{latfund3}
Suppose $U$ is a rational subspace of $V$ and $\Lambda_U =  U \cap \Lambda.$ (Note that $\Lambda_U$ is a lattice of $U.$) Then any basis of $\Lambda_U$ can be extended to a basis of $\Lambda.$
\end{lem}

The next lemma discusses the lattice of the projection of a subspace and the lattice of the kernel of the projection.
\begin{lem}\label{latfund4}
Suppose $\Gamma$ is a sublattice of $\Lambda$ and $U$ is the subspace of $V$ spanned by $\Gamma.$ Then $\pi^{(D-k)}(\Gamma)$ is a lattice of $\pi^{(D-k)}(U)$ and $\Lambda^k \cap \Gamma$ is a lattice of $V^k \cap U.$ 

Let $r = \dim(\pi^{(D-k)}(U))$ and $t = \dim(V^k \cap U).$ Then $r + t = \dim(U).$

Furthermore, suppose $(\bf^1, \dots, \bf^t)$ is a basis of the lattice $\Lambda^k \cap \Gamma$ and $(\bg^1, \dots, \bg^r)$ is a set of points in $\Gamma.$ Then $(\bg^1, \dots, \bg^r, \bf^1, \dots, \bf^t)$ is a basis of $\Gamma$ if and only if $(\pi^{(D-k)}(\bg^1), \pi^{(D-k)}(\bg^2), \dots, \pi^{(D-k)}(\bg^r))$ is a basis of the lattice $\pi^{(D-k)}(\Gamma).$
\end{lem}

\begin{proof}

$V^k \cap U$ is the kernel of the projection $\pi^{(D-k)}: U \to \pi^{(D-k)}(U).$ Thus, we can decompose the lattice $\Gamma$ as $\pi^{(D-k)}(\Gamma) \oplus (\Lambda^k \cap \Gamma).$
\end{proof}

\begin{lem}\label{charCentP}
Suppose $P$ is a $d$-dimensional central rational polytope in $V.$ 
\begin{ilist}
\itm $\Lambda_P = \aff(P) \cap \Lambda$ is a lattice of $\aff(P) = \lin(P).$
\itm $\pi^{(D-k)}(\Lambda_P)$ is a lattice of $\pi^{(D-k)}(\aff(P)).$
\itm $\Lambda^k \cap \Lambda_P$ is a lattice of $V^k \cap \aff(P).$
\itm $\dim(\pi^{(D-k)}(\aff(P))) + \dim(V^k \cap \aff(P)) = d.$
\end{ilist}
\end{lem}

\begin{proof}
Since $P$ is central and rational, we have that $\aff(P) = \lin(P)$ is rational. Thus, (i) follows from Lemma \ref{latfund1}, and (ii), (iii) and (iv) follow from Lemma \ref{latfund4}.
\end{proof}

\subsection{On integrality and general position}

We say an $\ell$-dimensional subspace $U$ of $V$ is the {\it row space} of an $\ell \times D$ matrix if the rows of the matrix are the coordinates of a basis of $U$ with respect to $\e.$

Recall a {\it unimodular matrix} is a square integer matrix with determinant $+1$ or $-1.$

The following two lemmas can be checked directly from the definition.
\begin{lem}\label{charInt}
Let $U$ be an $\ell$-dimensional subspace of $V.$ Then the following are equivalent:
\begin{ilist}
\itm $U$ is integral.
\itm $U$ is the row space of a matrix of the form $(I \ J)$, where $I$ is the $\ell \times \ell$ identity matrix and $J$ is an $\ell \times (D - \ell)$ integer matrix.
\itm $U$ is the row space of a matrix of the form $(K \ J)$, where $K$ is an $\ell \times \ell$ unimodular matrix and $J$ is an $\ell \times (D - \ell)$ matrix of integer entries. 
\end{ilist}

An affine space $W$ is integral if and only if $W$ contains a lattice point and $\lin(W)$ is integral if and only if $\lin(W)$ is integral and $W = \lin(W) + \z$ for some (equivalently all) $\z \in W \cap \Lambda.$
\end{lem}

\begin{lem}\label{charGen}
Let $U$ be an $\ell$-dimensional subspace of $V.$ Then the following are equivalent:
\begin{ilist}
\itm $U$ is in general position.
\itm $U$ is the row space of a matrix of the form $(I \ J)$, where $I$ is the $\ell \times \ell$ identity matrix and $J$ is an arbitrary $\ell \times (D - \ell)$ matrix. 
\itm $U$ is the row space of a matrix of the form $(K \ J)$, where where $K$ is an invertible matrix and $J$ is an arbitrary $\ell \times (D - \ell)$ matrix. 
\itm If $U$ is the row space of a matrix of the form $(K \ J)$, where $J$ is any arbitrary $\ell \times (D - \ell)$ matrix, then $K$ is an invertible matrix. 
\end{ilist}

An affine space $W$ is in general position if and only if $\lin(W)$ is in general position.
\end{lem}

\begin{cor}\label{propIntGen}
Suppose $0 \le \ell \le m \le D$ and $W$ is an $\ell$-dimensional affine space of $V.$
\begin{ilist} 
\itm If $W$ is integral, then $\lin(W)$ is rational and $W$ is in general position.
\itm If $W$ is integral, then the projection $\pi^{(D-m)}(W)$ is $\ell$-dimensional and integral. In particular, $\pi^{(D-\ell)}$ induces a bijection between $W \cap \Lambda$ and $\Lambda_\ell.$
\itm If $W$ is in general position, then the projection $\pi^{(D-m)}(W)$ is $\ell$-dimensional and in general position. In particular, $\pi^{(D-\ell)}$ induces a bijection between $W$ and $V_\ell.$
\itm If $W$ is in general position, for any $\ell+1$ affinely independent points $\x_1, \dots,$ $\x_{\ell+1} \in W,$ we have that the matrix $\left(\begin{array}{cc} 1 & \pi^{(D-\ell)}(\x_1) \\ 1 & \pi^{(D-\ell)}(\x_2) \\ \vdots & \vdots \\ 1 & \pi^{(D-\ell)}(\x_{\ell+1})\end{array}\right)$ is invertible.
\end{ilist}
\end{cor}

\begin{proof}
(i), (ii) and (iii) directly follows from Lemma \ref{charInt} and Lemma \ref{charGen}. One checks that
$$\det\left(\begin{array}{cc} 1 & \pi^{(D-\ell)}(\x_1) \\ 1 & \pi^{(D-\ell)}(\x_2) \\ \vdots & \vdots \\ 1 & \pi^{(D-\ell)}(\x_{\ell+1})\end{array}\right) = (-1)^\ell \det\left(\begin{array}{c} \pi^{(D-\ell)}(\x_1) -\pi^{(D-\ell)}(\x_{\ell+1})\\  \pi^{(D-\ell)}(\x_2)-\pi^{(D-\ell)}(\x_{\ell+1}) \\ \vdots  \\  \pi^{(D-\ell)}(\x_{\ell}) - \pi^{(D-\ell)}(\x_{\ell+1})\end{array}\right),$$
and $\lin(W)$ is the row space of the matrix $(\x_i - \x_{\ell+1})_{1\le i \le \ell},$ whose first $\ell$ columns are exactly the latter matrix in the above formula. Then (iv) follows from Lemma \ref{charGen}/(iv).
\end{proof}

Given (i) of Corollary \ref{propIntGen}, we can apply all the results we obtain for affine spaces (or polytopes) in general position to affine spaces (or polytopes) that are integral. This is an important fact which we use implicitly throughout the paper.

\begin{lem}\label{charRank}
Suppose $P \subset V$ is a $d$-dimensional rational polytope. 
\begin{ilist}
\itm If one of the $k$-faces of $P$ is affinely integral, then $\pi^{(D-k)}(\Lambda_P) = \Lambda_k.$
\itm If $P$ is central and one of the $k$-faces of $P$ is in affinely general position, then $\pi^{(D-k)}(\Lambda_P)$ is of rank $k$ and $\Lambda^k \cap \Lambda_P$ is of rank $d-k.$
\end{ilist}
\end{lem}

\begin{proof}
(i) follows immediately from Corollary \ref{propIntGen}/(ii).

(ii) follows from Corollary \ref{propIntGen}/(iii) and Lemma \ref{charCentP}.
\end{proof}

\begin{lem}\label{preDilateInt}
Suppose $P$ is a $k$-integral polytope. Then $m P$ is a $k$-integral polytope, for any positive integer $m.$
\end{lem}

\begin{proof}
This can be checked using Lemma \ref{charInt}.
\end{proof}

\subsection{Algebra of polyhedra}
Both the number of lattice points and the volume of a polytope are {\it valuations;} that is, they satisfy the inclusion-exclusion property. Since the theory of valuations is based on the theory of the algebra of polyhedra, we will recall the basic definitions in the area. Please refer to \cite{BarviPom} for details.

\begin{defn}
Let $A \subset V$ be a set. The {\it indicator function}, or just {\it indicator}, of $A$ is the function $[A]: V \to \R$ defined by
$$[A](x) = \begin{cases}
1 & \mbox{if $x \in A$}, \\
0 & \mbox{if $x \not\in A$}.
\end{cases}$$
The {\it algebra of polyhedra} $\cP(V)$ is the vector space over $\Q$ spanned by the indicator functions $[P]$ of all polyhedra $P \subset V.$
The {\it algebra of polytopes} $\cP_b(V)$ is the vector space over $\Q$ spanned by the indicator functions $[P]$ of all polytopes $P \subset V.$
\end{defn}

\begin{defn}
A linear transformation $\Phi: \cP(V) \to V',$ or $\Phi: \cP_b(V) \to V',$ where $V'$ is a vector space, is called a {\it valuation}.
\end{defn}

It is clear that one can extend the number of lattice points of a polytope to a valuation of $\cP_b(V).$ However, we have to be more careful when we extend the volume function to a valuation, because if we use different measures for different polytopes, the linearity will fail. Thus, we fix a lattice in the following definition. For convenience, we also fix a notation related to the right side of \eqref{VolFor}.

\begin{defn}
Let $\Gamma$ be a sublattice of $\Lambda$ and $U$ the subspace of $V$ spanned by $\Gamma.$ For any polytope $P \subset U,$ we define the {\it $k$th S-volume of $P$ (with respect to $\Gamma$)} to be
$$\svol_{\Gamma}^k(P) = \sum_{\y \in \pi^{(D-k)}(\Gamma)} \vol_{\Lambda^k \cap \Gamma}(\pi_{D-k}(\y,P)).$$ 
Recall $\vol_{\Gamma}(P)$ is defined by \eqref{defVol}.
One can extend both $\vol_{\Gamma}$ and $\svol_{\Gamma}^k$ to valuations of $\cP_b(U).$

For convenience, for any $[S] \in \cP_b(U),$ we omit the bracket and write $\vol_{\Gamma}(S)$ and $\svol_{\Gamma}^k(S)$ instead of $\vol_{\Gamma}([S])$ and $\svol_{\Gamma}^k([S]).$
\end{defn}

Using the definition of $\svol_\Gamma^k$ and Lemma \ref{cent}, we can rephrase Theorem \ref{mainVol} as the following proposition.
\begin{prop}\label{mainVol1}
Suppose $0 < k < d \le D$ and $P \subset V$ is a $d$-dimensional central $(k-1)$-integral polytope in $k$-general position. Let $\Lambda_P := \aff(P) \cap \Lambda.$ Then the volume of $P$ is given by
\begin{equation}\label{VolFor1}
\vol_{\Lambda_P}(P) = \svol_{\Lambda_P}^k(P).
\end{equation}
\end{prop}

In the next three sections, we will prove a special case of Proposition \ref{mainVol1}: when $k=1$ and $P$ is a simplex. The material in these sections is not used in the subsequent sections, so the reader who is willing to take this case of Proposition \ref{mainVol1} on faith may skip to Section \ref{propertySP}.

All of the results in Sections \ref{preservation} and \ref{reduction} and some of the results in Section \ref{special} are stated for general $k,$ since the statements and proofs are no harder than the $k=1$ case.

\section{Affine transformations that preserve certain properties}\label{preservation}
In this section, we will discuss conditions for affine transformations that preserve certain properties of polytopes: integrality, general position, volumes and S-volumes. This serves as a preparation for the reduction we will discuss in the next section.

Suppose $\Gamma_i$ is a sublattice of $\Lambda$ which has a basis $\e_i = (\be_i^{1}, \dots, \be_i^{D_i}),$ and $U_i$ the subspace of $V$ spanned by $\Gamma_i$, for $i=1,2.$ For any affine transformation $\phi: U_1 \to U_2$, we can associate a vector $\alpha^\phi_{\e_1, \e_2} \in \R^{D_2}$ and a $D_1 \times D_2$ matrix $M^\phi_{\e_1, \e_2}$ to $\phi$ such that $\phi$ sends a point in $U_1$ with coordinates $\x$ with respect to the basis $\e_1$ to the point in $U_2$ with coordinates $\alpha^\phi_{\e_1, \e_2} + \x M^\phi_{\e_1, \e_2}$ with respect to the basis $\e_2.$ One verifies that 
$$M^\phi_{\e_1, \e_2} = (m_{i,j}), \mbox{if }  \phi(\be_1^{i})  - \phi(0) = \sum_{j=1}^{D_2} m_{i,j} \be_2^{j},$$
and $\alpha^\phi_{\e_1, \e_2}$ are the coordinates of $\phi(0)$ with respect to $\e_2.$

If the map is from $V$ to $V,$ unless otherwise noted we will assume that the lattice is $\Lambda$ with the fixed basis $\e = (\be^1, \dots, \be^D)$ and omit the subscript ``$\e_1, \e_2$''. Similarly, we can associate a matrix to any linear transformation.

\subsection{Conditions for preserving integrality and general position}
We say an affine map $\phi$ {\it preserves integrality} of an integral affine space $W$ if $\phi(W)$ is integral. Similarly, we say an affine map $\phi$ {\it preserves generality} of an affine space $W$ in general position, if $\phi(W)$ is in general position.

\begin{lem}\label{preIntGen}
Suppose $\phi: V \to V$ is an affine map with associated matrix
$$M = \left(\begin{array}{cc} A_{\ell \times \ell} & C_{\ell \times (D-\ell)} \\ O_{(D-\ell) \times \ell} & B_{(D-\ell) \times (D-\ell)} \end{array} \right)$$ and vector $\alpha,$ where $O$ is the zero matrix and $B$ is invertible. 
\begin{ilist}
\itm Let $W$ be an $\ell$-dimensional integral affine space. If $A$ is unimodular and $\phi$ sends any lattice point in $W$ to a lattice point, then $\phi$ preserves integrality of $W.$
\itm $\phi$ preserves integrality of any $\ell$-dimensional integral affine spaces of $V$ if $A$ is unimodular, and both $B$ and $C$ are integer matrices, and $\alpha$ is an integer vector.
\itm $\phi$ preserves generality of any  $\ell$-dimensional affine spaces of $V$ in general position if $A$ is invertible.
\end{ilist}
\end{lem}

\begin{proof}
(ii) follows from (i), and (iii) can be check directly using Lemma \ref{charGen}/(iii). Now we prove (i). Let $U = \lin(W),$ then by Lemma \ref{charInt}, $U$ is integral and $W = U + \z,$ for some $\z \in W \cap \Lambda.$ By Lemma \ref{charInt}/(ii), $U$ is the row space of an integer matrix of the form $(I \ J),$ where $I$ is the $\ell \times \ell$ identity matrix. One checks that $\phi(W) = U' + \phi(\z),$ where $U'$ is the row space of $(I \ J)M = (A \ C+JB)$ and $\phi(\z) \in \phi(W) \cap \Lambda.$ Thus, it is enough to show that $U'$ is integral. Since $A$ is unimodular, by Lemma \ref{charInt}/(iii), we only need to verify that $C+JB$ is an integer matrix, or equivalently each row in $(I \ J)M$ is an integer vector. Thus, we need to show that for any row vector $\v$ in $(I \ J),$ $\v M$ is an integer vector. However, $\v$ is an integer vector and $W$ contains lattice points, so $\v$ can be written as the difference of two integer vectors $\u_1$ and $\u_2,$ each of which is the coordinates of a lattice point in $W.$ Because $\phi$ sends lattice points in $W$ to lattice points, we have that $\alpha + \u_i M$ is an integer vector, for each $i.$ Therefore, $\v M$ is an integer vector as well. 
\end{proof}

The following corollary follows immediately.

\begin{cor}\label{preIntGen1}
Suppose $\phi: V \to V$ is an invertible affine map with an associated matrix
$$M = \left(\begin{array}{cc} A_{k \times k} & C_{k \times (D-k)} \\ O_{(D-k) \times k} & B_{(D-k) \times (D-k)} \end{array} \right)$$ and vector $\alpha$ where $O$ is the zero matrix. (We must have that $A$ and $B$ are invertible.) 
\begin{ilist}
\itm Suppose $W$ is an affine space of $V$ satisfying that $\phi$ sends any lattice point in $W$ to a lattice point. If $A$ is an upper triangular and unimodular matrix, then $\phi$ preserves integrality of any integral affine space of $W$ of dimension no more than $k$.
\itm If $A$ is an upper triangular and unimodular matrix,  both $B$ and $C$ are integer matrices, and $\alpha$ is an integer vector, then $\phi$ preserves integrality of any integral affine space of $V$ of dimension no more than $k$.
\itm If $A$ is upper triangular, then $\phi$ preserves generality of any affine space of $V$ of dimension no more than $k$ that is in general position.
\end{ilist}

\end{cor}

\subsection{Conditions for preserving volumes and S-volumes}
We will discuss situations when $\vol_{\Gamma}(P)$ and $\svol_{\Gamma}^k(P)$ are invariant under affine transformation.  

\begin{lem}\label{invar}
Let $s$ be a positive integer. For $i=1,2,$ let $\Gamma_i$ be a sublattice of $\Lambda$ of rank $s$ and $U_i$ the subspace of $V$ spanned by $\Gamma_i$. Suppose $\e_i = (\be_i^{1}, \dots, \be_i^{s})$ is a basis of $\Gamma_i,$ for $i=1,2.$ Let  $\phi: U_1 \to U_2$ be an invertible affine map with the associated matrix $M_{\e_1,\e_2}$.
\begin{ilist}
\itm If $\det(M_{\e_1, \e_2}) =1,$ then $\vol_{\Gamma_1} = \vol_{\Gamma_2} \circ \phi.$
\itm Suppose the following four conditions conditions are satisfied:
\begin{enumerate}
 \item $\pi^{(D-k)}(\x_1) = \pi^{(D-k)}(\x_2)$ if and only if $\pi^{(D-k)}(\phi(\x_1))= \pi^{(D-k)}(\phi(\x_2)).$
\item $\pi^{(D-k)}(\x) \in \pi^{(D-k)}(\Gamma_1)$ if and only if $\pi^{(D-k)}(\phi(\x)) \in \pi^{(D-k)}(\Gamma_2).$
\item For some number $r$ and for $i=1,2,$ we have that $(\be_i^{r+1}, \be_i^{r+2}, \dots, \be_i^{s})$ is a basis of the lattice $\Lambda^k \cap \Gamma_i.$ (Thus, it spans $V^k \cap U_i$ by Lemma \ref{latfund2}.)
\item The determinant of the lower right $(s-r) \times (s-r)$ submatrix of $M_{\e_1, \e_2}$ is equal to $1.$
\end{enumerate}
Then $\svol_{\Gamma_1}^k = \svol_{\Gamma_2}^k \circ \phi.$
\end{ilist} 
\end{lem}

\begin{proof}
(i) follows directly from the definition of volume and linear algebra.

By (1), one sees that $\phi$ induces a bijection between $\pi^{(D-k)}(U_1)$ and $\pi^{(D-k)}(U_2).$ We denote this map by $\phi'.$ We then have that for any set $S \in U_1$ and $\y \in \pi^{(D-k)}(U_1),$
$$\pi_{D-k}(\phi'(\y), {\phi}(S)) =  {\phi} (\pi_{D-k}(\y, S)).$$
Note that shifting a set of points by a constant vector does not change its volume. Thus,
\begin{equation}\label{shift1}
\vol_{\Lambda^k \cap \Gamma_2}( {\phi} (\pi_{D-k}(\y, S))) = \vol_{\Lambda^k \cap \Gamma_2}( \{  \x - \phi'(\y) \ | \ \x \in {\phi} (\pi_{D-k}(\y, S))  \}),
\end{equation}
and
\begin{equation}\label{shift2}
\vol_{\Lambda^k \cap \Gamma_1} (\pi_{D-k}(\y, S)) = \vol_{\Lambda^k \cap \Gamma_1} (\{ \x - \y \ | \ \x \in \pi_{D-k}(\y, S) \}).
\end{equation}
We will show the right sides of the above two equations are equal.
We define an affine map $\phi_\y: V^k \cap U_1 \to V^k \cap U_2$ by mapping $\z$ to $\phi(\y + \z) - \phi'(\y).$ Using the bases $(\be_1^{r+1}, \be_1^{r+2}, \dots, \be_1^{s})$ and $(\be_2^{r+1}, \be_2^{r+2}, \dots, \be_2^{s}),$ one can check the matrix associated to $\phi_\y$ is exactly the lower right $(s-r) \times (s-r)$ submatrix of $M_{\e^1, \e^2}.$ Hence, by (4), we have that 
$\vol_{\Lambda^k \cap \Gamma_1} = \vol_{\Lambda^k \cap \Gamma_2} \circ \phi_\y.$
We verify that $$\phi_\y(\{ \x - \y \ | \ \x \in \pi_{D-k}(\y, S) \}) = \{  \x - \phi'(\y) \ | \ \x \in {\phi} (\pi_{D-k}(\y, S))  \}.$$
Therefore, the right sides of equations \eqref{shift1} and \eqref{shift2} are equal. So the left sides are equal, and we conclude that for any set $S \in U_1$ and $\y \in \pi^{(D-k)}(U_1),$
\begin{eqnarray*}
\vol_{\Lambda^k \cap \Gamma_2} (\pi_{D-k}(\phi'(\y), {\phi}(S)))  &=& \vol_{\Lambda^k \cap \Gamma_2}( {\phi} (\pi_{D-k}(\y, S))) \\
&=& \vol_{\Lambda^k \cap \Gamma_1} (\pi_{D-k}(\y, S)).
\end{eqnarray*}
Furthermore, by (2), $\phi'$ also gives a bijection between $\pi^{(D-k)}(\Gamma_1)$ and $\pi^{(D-k)}(\Gamma_2).$ 
Thus, we have $\svol_{\Gamma_1}^k = \svol_{\Gamma_2}^k \circ \phi.$
\end{proof}

\begin{cor}\label{invarVol1}
Let $s$ be a positive integer. For $i=1,2,$ let $\Gamma_i$ be a sublattice of $\Lambda$ of rank $s$ and $U_i$ the subspace of $V$ spanned by $\Gamma_i.$ Suppose $\e_i = (\be_i^{1}, \dots, \be_i^{s})$ is a basis of $\Gamma_i,$ and for some number $r$ we have that $(\be_i^{r+1}, \be_i^{r+2}, \dots, \be_i^{s})$ is a basis of the lattice $\Lambda^k \cap \Gamma_i.$ Suppose ${\phi}: U_1 \to U_2$ is an affine map with an associated matrix
$$M^\phi_{\e_1, \e_2} = \left(\begin{array}{cc} A & C \\ O & B \end{array} \right)$$ and vector $\alpha_{\e^1, \e^2} = (\alpha_1, \alpha_2)$ satisfying
\begin{ilist}
\itm $A$ is an $r \times r$ unimodular matrix.
\itm $\alpha_1$ is an $r$-dimensional vector of integer entries.
\itm $B$ is an $(s-r) \times (s-r)$ matrix with $\det(B) = 1.$
\itm $O$ is the zero matrix.
\end{ilist}
Then $\vol_{\Gamma_1} = \vol_{\Gamma_2} \circ \phi$ and $\svol_{\Gamma_1}^k = \svol_{\Gamma_2}^k \circ \phi.$

\end{cor}

\begin{proof}
It is clear that $\det(M_{\e_1, \e_2}) =1.$ Thus, $\vol_{\Gamma_1} = \vol_{\Gamma_2} \circ \phi.$ To show that $\svol_{\Gamma_1}^k = \svol_{\Gamma_2}^k \circ \phi$, we need to verify (1)--(4) under (ii) of Lemma \ref{invar} holds. (3) and (4) are given in the conditions. 

By Lemma \ref{latfund4}, $\f_i =(\pi^{(D-k)}(\be_i^1), \dots, \pi^{(D-k)}(\be_i^r))$ is a basis of $\pi^{(D-k)}(\Gamma_i),$ for $i=1,2.$ One checks that $\phi$ deduces an affine map from $\pi^{(D-k)}(U_1)$ to $\pi^{(D-k)}(U_2)$ by mapping $\pi{(D-k)}(\x)$ to $\pi{(D-k)}(\phi(\x)).$ The matrix and vector associated to this map with respect to the bases $\f_1$ and $\f_2$ are $A$ and $\alpha_1,$ respectively. Therefore, (1) and (2) of Lemma \ref{invar} are satisfied.
\end{proof}


\section{Reduction}\label{reduction}

In this section, we will prove the following proposition which reduce the problem of proving Proposition \ref{mainVol1} to full-dimensional polytopes in fully general position.
\begin{prop}\label{red2full}
Suppose $0 < k < d \le D$ and $P \subset V$ is a $d$-dimensional central $(k-1)$-integral polytope in $k$-general position. There exists an invertible affine map $\phi: V \to V$ such that $Q := \phi(P)$ is a full-dimensional polytope in $V_d$ that is both $(k-1)$-integral polytope and in fully general position. Also, proving Proposition \ref{mainVol1} for $P$ is equivalent to proving it for $Q.$
\end{prop}

Throughout this section, we will assume that $0 < k < d \le D$ and $P \subset V$ is a $d$-dimensional central $(k-1)$-integral polytope in $k$-general position. We also fix
$$U_P := \aff(P),$$ which is a subspace of $V$ because $P$ is central.

\subsection{Reduction to the full-dimensional case}
\begin{lem}\label{reddim}
There exists an invertible linear transformation $\phi: V \to V$ such that $Q := \phi(P)$ is a full-dimensional polytope in $V_d$ that is still $(k-1)$-integral and in $k$-general position. Also, proving Proposition \ref{mainVol1} for $P$ is equivalent to proving it for $Q.$
\end{lem}

We first prove two preliminary lemmas.
\begin{lem}\label{charLamP}
There exist integers $h_1, \dots, h_k$ such that
$(\be^1 + h_1 \be^k ,\dots, \be^{k-1} + h_{k-1} \be^{k}, h_k \be^k)$ is a basis of $\pi^{(D-k)}(\Lambda_P).$
\end{lem}
\begin{proof}
By Lemma \ref{charRank}, we have that $\pi^{(D-(k-1))}(\Lambda_P) = \Lambda_{k-1},$ and $\pi^{(D-k)}(\Lambda_P)$ is of rank $k$. An obvious basis for $\pi^{(D-(k-1))}(\Lambda_P)$ is $(\be^1, \be^2, \dots, \be^{k-1}).$ Note that $\pi^{(D-(k-1))}(\Lambda_P)$ can be obtained by dropping the last coordinate, i.e., the coordinate corresponding to $\be^k,$ of the lattice $\pi^{(D-k)}(\Lambda_P).$ Hence, our lemma follows from Lemma \ref{latfund4}.
\end{proof}

\begin{lem}\label{newbs1}
There exists a basis $\f = (\bf^1, \bf^2, \dots, \bf^D)$ of $V$ satisfying
\begin{ilist}
\itm $\pi^{(D-k)}(\bf^i) = \be^i + h_i \be^k,$ for $i =1, \dots, k-1,$ and $\pi^{(D-k)}(\bf^k) = h_k \be^k,$ for some integers $h_1, \dots, h_k.$ 
\itm $(\bf^{k+1}, \dots, \bf^d)$ is a basis of the lattice $\Lambda^k \cap \Lambda_P$
\itm $(\bf^1, \dots, \bf^k, \bf^{k+1}, \dots \bf^d)$ is a basis of the lattice $\Lambda_P.$
\itm $(\bf^{k+1}, \dots, \bf^d, \bf^{d+1}, \dots,\bf^D)$ is a basis of the lattice $\Lambda^k.$
\end{ilist}
\end{lem}

\begin{proof}
By Lemma \ref{charLamP}, $\pi^{(D-k)}(\Lambda_P)$ has a basis $$(\be_1 + h_1 \be^k ,\dots, \be_{k-1} + h_{k-1} \be^{k-1}, h_k \be^k),$$ for some integers $h_1, \dots, h_k.$
Thus, using Lemma \ref{latfund4} again, we can choose $\bf^1, \dots, \bf^{k},$ $\bf^{k+1}, \dots, \bf^d$ satisfying (i), (ii) and (iii).

Furthermore, by Lemma \ref{latfund3}, $(\bf^{k+1}, \dots, \bf^d)$ can be extended to a basis, say $(\bf^{k+1}, \dots, \bf^d, \bf^{d+1}, \dots, \bf^D)$ of $\Lambda^k.$
\end{proof}

\begin{proof}[Proof of Lemma \ref{reddim}]
Let $\f = (\bf^1, \dots, \bf^D)$ be a basis as described in Lemma \ref{newbs1}, and define $\phi: V \to V$ to be the linear transformation that maps $\bf^i$ to $\be^i$ for each $i.$ Both $\f$ and $\e$ are bases of $V.$ Thus, $\phi$ is invertible. Let $Q = \phi(P).$ Clearly, $Q$ is a $d$-dimensional polytope in $V_d.$ 

One checks that the associated matrix to $\phi$ is $M^\phi = M^\phi_{\e, \e} = F^{-1},$ where $F$ is the $D \times D$ matrix whose $i$th row are the coordinates of $\bf^i$ with respect to $\e.$ By Lemma \ref{newbs1}, we know that $F$ is a matrix satisfying the following three properties:
\begin{ilist}
\itm The lower left $(D-k) \times k$ submatrix of $F$ is a zero matrix.
\itm The upper left $k \times k$ submatrix of $F$ is an upper triangular matrix.
\itm The upper left $(k-1) \times (k-1)$ submatrix of $F$ is an upper triangular unimodular matrix.
\end{ilist}
It is easy to verify that $M^\phi = F^{-1}$ also has these three properties. Recall that $U_P = \aff(P)$ is a subspace of $V,$ and $\Lambda_P = \aff(P) \cap \Lambda.$ Hence, by the construction of $\phi,$ any lattice point in $U_P$ is mapped to a lattice point by $\phi.$  Therefore, by using (i) and (iii) of Corollary \ref{preIntGen1}, we conclude that $Q$ is a $(k-1)$-integral polytope in $k$-general position.

The map $\phi$ induces a linear map $\tilde{\phi}: U_P \to V_d.$ It is clear that the associated matrix $M^{\tilde{\phi}}_{\tilde{\f}, \tilde{\e}}$ to $\tilde\phi$ with respect to the bases $\tilde{\f} = (\bf^1, \dots, \bf^d)$ and $\tilde{\e} = (\be^1, \dots, \be^d)$  is the identity matrix.
Note that by (ii) of Lemma \ref{newbs1}, $(\bf^{k+1}, \dots, \bf^{d})$ is a basis of $\Lambda^k \cap \Lambda_P,$ and it is obvious that $(\be^{k+1}, \dots, \be^{d})$ is a basis of $\Lambda^k \cap \Lambda_d.$ Thus, by Corollary \ref{invarVol1}, $\vol_{\Lambda_P}(P) = \vol_{\Lambda_d}(Q)$ and $\svol_{\Lambda_P}^k(P) = \svol_{\Lambda_d}^k(Q).$ However, $Q$ is full-dimensional in $V_d.$ Thus, $\Lambda_Q = V_d \cap \Lambda = \Lambda_d.$

Thus, we conclude that proving Proposition \ref{mainVol1} for $P$ is reduced to proving it for $Q = \phi(P).$ Therefore, Lemma \ref{reddim} follows the fact that $Q$ is full-dimensional in $V_d.$
\end{proof}

\subsection{Reduction to the fully general position case}
Because of Lemma \ref{reddim}, we will assume in this subsection that $P$ is a full-dimensional polytope in $V,$ or in other words that $d = D.$ 

\begin{lem}\label{redgen}
Suppose $0 < k \le \ell < d.$ Let $P$ be a $d$-dimensional $(k-1)$-integral polytope in $\ell$-general position in $V = V_d$. Then there exists an invertible linear transformation $\phi: V \to V$ such that $Q := \phi(P)$ is a $d$-dimensional $(k-1)$-integral polytope in $(\ell+1)$-general position, and $\vol_{\Lambda_P}(P) = \vol_{\Lambda_Q}(Q)$ and $\svol_{\Lambda_P}^k(P) = \svol_{\Lambda_Q}^k(Q).$
\end{lem}

\begin{lem}\label{findw}
For any finite set of nonzero vectors $\{ \v_i \}$ in $\R^{m},$ there exists an integer vector $\w \in \Z^m$ with the first coordinate $1$ such that $\v_i \cdot \w \neq 0$ for each $i.$
\end{lem}
\begin{proof}
For each $i,$ let $H_i$ be the null space of $\v_i.$ 
Let $H_0$ be the set of points in $\R^m$ with the first coordinate equal to $1$. The set $\H := \cup_i (H_i \cap H_0)$ is a finite union of $(m-2)$-dimensional affine spaces inside $H_0.$ One sees easily that the complement (with respect to $H_0$) of $\H$ contains lattice points. 
We can let $\w$ be any such lattice point.
\end{proof}

\begin{proof}[Proof of Lemma \ref{redgen}]
Let $\{ U_i \}$ be the finite set of subspaces of $V$ that are translations of affine hull of $(\ell+1)$-faces of $P.$ Since $P$ is $\ell$-general position, each $U_i$ has an $\ell$-dimensional subspace that is in general position. Therefore, using Lemmas \ref{charGen}, we see that $U_i$ is the row space of a matrix of the form
$$\left(\begin{array}{cc} I &  J_i \\ {\mathbf o} & \v_i \end{array} \right),$$
where $I$ is the $\ell \times \ell$ identity matrix and $J_i$ is an $\ell \times (D - \ell)$ matrix, $\mathbf o$ is the $\ell$-dimensional zero vector, and $\v_i$ is a $(d-\ell)$-dimensional nonzero vector.

By Lemma \ref{findw}, there exists $\w \in \Z^{d-\ell}$ such that the first coordinate of $\w$ is $1$ and $\v_i \cdot \w \neq 0$ for each $i.$ Let $M$ be the $D \times D$ matrix obtained from the identity matrix by replacing the $(\ell+1)$th column with $({\mathbf o}, \w)^T.$ We consider the linear transformation $\phi: V \to V$ with associated matrix $M^\phi = M.$ The determinant of $M^\phi$ is $1,$ thus $\phi$ is invertible and $Q$ is $d$-dimensional.

It follows from Corollary \ref{preIntGen1}/(ii)(iii) that Q is a $(k-1)$-integral polytope in $\ell$-general position. Since both $P$ and $Q$ are full-dimensional, we have that $\Lambda_P = \Lambda_Q = \Lambda.$ Then, by Corollary \ref{invarVol1},  we have that $\vol_{\Lambda_P}(P) = \vol_{\Lambda_Q}(Q)$ and $\svol_{\Lambda_P}^k(P) = \svol_{\Lambda_Q}^k(Q).$

It is left to show that any $(\ell+1)$ face $F$ of $P$ is in general position. Suppose $F$ is an $(\ell+1)$ face of $P.$ Then $\lin(F) = U_i$ for some $i.$ One checks that $\lin(\phi(F)) = \phi(U_i).$  However, $\phi(U_i)$ is the row space of a matrix of the form
$$\left(\begin{array}{cc} I &  J_i \\ {\mathbf o} & \v_i \end{array} \right) M^{\phi} = \left(\begin{array}{cc} I &  J_i' \\ {\mathbf o} & \v_i' \end{array} \right),$$
where the first coordinate of $\v_i'$ is $\v_i \cdot \w \neq 0.$ Thus, by Lemma \ref{charGen}, we conclude that $\phi(U_i)$ is in general position, and thus so is $F.$
\end{proof}

\begin{proof}[Proof of Proposition \ref{red2full}]
The proposition follows from Lemma \ref{reddim} and recursively applying Lemma \ref{redgen}.
\end{proof}

\section{A special case of Proposition \ref{mainVol1}}\label{special}

In this section, we prove the following proposition as a special case of Proposition \ref{mainVol1}.
\begin{prop}\label{k=1simplex}
Propostion \ref{mainVol1} is true when $k=1,$ and $P$ is a simplex.
\end{prop}

By Proposition \ref{red2full}, we can prove Proposition \ref{k=1simplex} with the assumption that $d=D$ and $P$ is in fully general position.

\subsection{Basic setup and notation}
We will use quite a few results from \cite{lattice-face, note-lattice-face} in the proof of Proposition \ref{k=1simplex}. We start by recalling relevant definitions. We remark that although in \cite{lattice-face, note-lattice-face} we worked in the space $\R^d$ with lattice $\Z^d,$ because of the canonical maps between $V$ and $\R^d$ that maps $\Lambda$ to $\Z^d,$ we can easily apply results in \cite{lattice-face, note-lattice-face} to our setting.

For simplicity, we use $\pi$ to denote the projection from $V_k$ to $V_{k-1}$ obtained by dropping the last coordinate, for any $k: 1\le k \le d.$ 
For any set $S \subset V_k$ and any point $\y \in V_{k-1},$ let $\rho_k(\y, S) = \pi^{-1}(\y) \cap S$ be the intersection of $S$ with the inverse image of $\y$ under $\pi.$ If $S$ is bounded, let $n_k(\y, S)$ be the points in $\rho_k(\y,S)$ with the smallest last coordinate. If $\rho_k(\y,S)$ is the empty set, i.e., $\y \not\in \pi(S),$ then let $n_k(\y,S)$ be the empty set as well.

\begin{defn}\label{dfbdry}
Let $P$ be a polytope in $V_k.$ Define 
$$\Omega_k(P) = P \setminus \bigcup_{\y \in V_{k-1}} n_k(\y, P)$$ to be the {\it nonnegative part} of $P.$

If $P$ is full dimensional in $V_k,$ we often omit the subscript $k$ and just write $\Omega(P).$


\end{defn}



For the rest of this section, we will assume $D=d \ge 2,$ and $P \in V=V_d$ is a $d$-dimensional simplex in general position with vertex set $\vert = \{v_1, v_2, \dots, v_{d+1}\},$ where the coordinates of $v_i$ are
$ \x_i = (x_{i,1},
x_{i,2}, \dots, x_{i,d}).$ 

For any $\sigma \in \S_d$ and $k: 1 \le k \le d,$ we define matrices
$X_\vert(\sigma, k)$ and $Y_\vert(\sigma, k)$ to be the matrices
$$X_\vert(\sigma, k) = \left(\begin{array}{ccccc}
1  & x_{\sigma(1), 1} & x_{\sigma(1),2} & \cdots & x_{\sigma(1), k} \\
1  & x_{\sigma(2), 1} & x_{\sigma(2),2} & \cdots & x_{\sigma(2), k} \\
\vdots  & \vdots & \vdots & \ddots & \vdots \\
1  & x_{\sigma(k), 1} & x_{\sigma(k),2} & \cdots & x_{\sigma(k), k} \\
1  & x_{d+1, 1} & x_{d+1,2} & \cdots & x_{d+1, k}
\end{array}\right)_{(k+1) \times (k+1)},$$

$$Y_\vert(\sigma, k) = \left(\begin{array}{ccccc}
1  & x_{\sigma(1), 1} & x_{\sigma(1),2} & \cdots & x_{\sigma(1), k-1} \\
1  & x_{\sigma(2), 1} & x_{\sigma(2),2} & \cdots & x_{\sigma(2), k-1} \\
\vdots  & \vdots & \vdots & \ddots & \vdots \\
1  & x_{\sigma(k), 1} & x_{\sigma(k),2} & \cdots & x_{\sigma(k),
k-1}
\end{array}\right)_{k \times k}.$$
We also define $z_\vert(\sigma, k)$ to be
$$z_\vert(\sigma,k) = \det(X_\vert(\sigma, k))/\det(Y_\vert(\sigma, k)).$$ Note that because $P$ is in fully general position and is a simplex, by Corollary \ref{propIntGen}/(iv), we have that $\det(X_\vert(\sigma, k))$ and $\det(Y_\vert(\sigma, k))$ are nonzero. Hence, $z_\vert(\sigma, k)$ is well-defined and is nonzero.

We often omit the subscript $\vert$ for $X_\vert(\sigma, k),$ $Y_\vert(\sigma,
k)$ and $z_\vert(\sigma,k)$ if there is no confusion.

Note that since $\Lambda_P = \Lambda,$ to prove Proposition \ref{k=1simplex}, we want to show that
$$\vol_{\Lambda}(P) = \svol^1_{\Lambda}(P).$$
However, since $P$ is a simplex, we have that
$$\vol_{\Lambda}(P) = \left|\frac{1}{d!}\det(X(\1, d)) \right|,$$ where $\1$ is the identity in $\S_d.$
Our goal is to describe $\svol^1_{\Lambda}(P)$ also in terms of $X(\sigma, k)$, $Y(\sigma,k)$ or $z(\sigma,k).$
We first have the following lemma.

\begin{lem}\
\begin{equation}\label{svol2Omega}
\svol^1_{\Lambda}(P) = \svol^1_{\Lambda}(\Omega(P)).
\end{equation}
\end{lem}

\begin{proof}
Since $\pi^{(d-1)}(\Omega(P)) \subset \pi^{(d-1)}(P),$ we can write
$$\svol^1_{\Lambda}(\Omega(P)) = \sum_{\y \in \Lambda_1 \cap \pi^{(d-1)}(P)} \vol_{\Lambda^1}(\pi_{d-1}(\y, \Omega(P))).$$
One sees that $\pi_{d-1}(\y, \Omega(P))) = \Omega_d(\pi_{d-1}(\y,P))$ which is obtained from $\pi_{d-1}(\y,P)$ by removing part of its boundary. Therefore, we must have that
$$\vol_{\Lambda^1}(\pi_{d-1}(\y, \Omega(P))) = \vol_{\Lambda^1}(\pi_{d-1}(\y, P)).$$
Hence, \eqref{svol2Omega} follows.
\end{proof}

Therefore, our problem becomes to find a way to describe $\svol^1_{\Lambda}(\Omega(P)).$

\subsection{Signed decomposition}
We use results in \cite{lattice-face} to find formulas for $\svol^1_{\Lambda}(\Omega(P)).$

\begin{lem}[Theorem 4.6 in \cite{lattice-face}]\label{decomp} 
Let $P \subset V = V_d$ be a $d$-dimensional simplex in fully general position with vertex set
$\vert = \{v_1, v_2, \dots, v_{d+1}\},$ where the coordinates of $v_i$ (with respect to $\e$)
are $ \x_i = (x_{i,1}, x_{i,2}, \dots, x_{i,d}).$ For any $\sigma
\in \S_d,$ and $k: 0 \le k \le d-1,$ let $v_{\sigma, k}$ be the
point with the first $k$ coordinates the same as $v_{d+1}$ and affinely
dependent with $v_{\sigma(1)}, v_{\sigma(2)}, \dots, v_{\sigma(k)},
v_{\sigma(k+1)}.$ (By Lemma 4.2/{\rm (vii)} in \cite{lattice-face}, we know that
there exists one and only one such point.) We also let $v_{\sigma,
d} = v_{d+1}.$ Then
\begin{equation}\label{fdecomp} 
[\Omega(P)] = \sum_{\sigma \in \S_d} \sgn(\sigma, P) [S_{\sigma}],
\end{equation}
 where
\begin{equation}\label{sign} \sgn(\sigma, P) = \sgn(\det(X(\sigma,
d))) \sgn\left(\prod_{i=1}^d z(\sigma, i)\right),
\end{equation}
and \begin{equation}\label{Ssigma} S_{\sigma}= \{ \s \in V \ | \
\pi^{(d-k)}(\s) \in \Omega(\pi^{(d-k)}(\conv(\{v_{\sigma,0}, \dots,
v_{\sigma,k}\}))) \forall 1 \le k \le d\}. \end{equation}

\end{lem}

We remark that multisets were used in \cite{lattice-face} to keep track of how many times each point is counted, we use indicator functions instead here. 

Since $\svol_{\Lambda}^1$ is a valuation, we reduce our problem to finding formulas for $\svol_{\Lambda}^1(S_\sigma).$ However, the structure of $S_\sigma$ is not quite good enough, so we look for affine transformations that preserves the S-volume but change $S_\sigma$ to a set whose S-volume is relatively easier to calculate.

In the proof of Proposition 5.2 of \cite{lattice-face}, an affine transformation $T_\sigma$ based on the vertex set of $P$ was defined by sending a point with coordinates $\x$ to a point with coordinates $\alpha_\sigma + \x M_\sigma.$ The construction of $\alpha_\sigma$ and $M_\sigma$ is not essential for our paper, so we won't state them here. Instead we will state properties of $T_\sigma$ that are useful for our purpose. We state below two simple facts of $T_\sigma$, which can be observed by checking the definitions directly, and by using Proposition 5.4 of \cite{lattice-face}.

\begin{lem}\label{Tsigma}
The affine transformation $T_\sigma$ defined in the proof of Proposition 5.2 of \cite{lattice-face} satisfies: 
\begin{ilist}
\itm The matrix $M_\sigma$ associated to $T_\sigma$ is upper-triangular with $1$'s on the diagonal. 
\itm The first coordinate of the constant vector $\alpha_\sigma$ associated to $T_\sigma$ is $-x_{\sigma(1),1}.$
\end{ilist}

Let $S_\sigma$ be defined as in Lemma \ref{decomp} and $\nHat{S}_{\sigma} = T_{\sigma}(S_{\sigma}).$ Then
\begin{equation}\label{Ssigmahat}
\nHat{S}_{\sigma} = \{ (s_1, s_2, \dots, s_d) \in V_d \ | \  s_k \in \Omega(\conv(0,
\frac{z(\sigma,k)}{z(\sigma,k-1)}s_{k-1} )), \forall 1 \le k \le d\},
\end{equation}
where by convention we let $z(\sigma, 0) = 1$ and $s_0 = 1.$

\end{lem}

\begin{rem}
Both Propositions 5.2 and 5.4 of \cite{lattice-face} were stated under the assumption that $P$ is a lattice-face polytope. The conclusion of Proposition 5.2 that $T_\sigma$ is a lattice-preserving map relies on this assumption, so is not always right if $P$ is just a polytope in fully general position. The two properties stated in Lemma \ref{Tsigma} can be viewed as a modified weaker version of Proposition 5.2 of \cite{lattice-face}. At the same time, Proposition 5.4 of \cite{lattice-face} does not rely on this assumption, so still holds under our current setting. 

\end{rem}

\begin{lem}\label{2hat}
Let $S_\sigma$ be defined as in Lemma \ref{decomp} and $\nHat{S}_{\sigma}$ as in Lemma \ref{Tsigma}. If $P$ is ($0$-)integral, then
for any $\sigma \in \S_d,$ we have that 
$$\svol_{\Lambda}^1(S_\sigma) = \svol_{\Lambda}^1(\nHat{S}_\sigma).$$
\end{lem}

\begin{proof}
Since $P$ is integral, we have that $-x_{\sigma(1),1}$ is always an integer. Hence, the lemma follows from Lemma \ref{invarVol1} and Lemma \ref{Tsigma}.
\end{proof}

\subsection{Power sums and formulas for $\svol_{\Lambda}^1(\nHat{S}_\sigma)$} We will describe $\svol_{\Lambda}^1(\nHat{S}_\sigma)$ in terms of $z(\sigma,k)$'s and power sums. Let's first review the definition of power sums and some of its properties.

Let $k$ be a positive integer. For any positive integer $x,$  we define
$$P_k(x) = \sum_{i=0}^x i^k.$$
It is well known that
\begin{eqnarray}
& &\mbox{$P_k(x)$ is a polynomial in $x$ of degree $k+1,$} \label{pdeg}\\
& &\mbox{the constant term of $P_k(x)$ is $0,$ i.e., $x$ is a factor
of $P_k(x),$} \label{pconst}\\
& &\mbox{the leading coefficient of $P_k(x)$ is $\frac{1}{k+1}.$}
\label{plead}
\end{eqnarray}
Because of \ref{pdeg}, we call $P_k(x)$ the {\it $k$th power sum polynomial.} Also, we can extend the domain of $P_k(x)$ from positive integers to real numbers. We have the following lemma from \cite{lattice-face}.

\begin{lem}[Lemma 5.8 in \cite{lattice-face}] For any positive integer $k,$
\begin{equation}\label{recp}
P_k(x) = (-1)^{k+1} P_k(-x-1).
\end{equation}
\end{lem}

Now we can describe the main result of this subsection.

\begin{lem}Let $\nHat{S}_{\sigma}$ be defined as in Lemma \ref{Tsigma}. If $P$ is ($0$-)integral, then
\begin{equation}\label{svolnHat}
\svol_{\Lambda}^1(\nHat{S}_\sigma) = \frac{1}{(d-1)!} \frac{\left|\prod_{j=1}^d z(\sigma,j)\right|}{z(\sigma,1)^{d}} P_{d-1}(z(\sigma,1)).
\end{equation}
\end{lem}

\begin{proof}
For convenience, for any two numbers $a \neq b,$ we denote by $[a, b] = \conv(a, b)$ the closed interval between $a$ and $b.$ Note that $\Omega([a,b]) = [a,b]\setminus \min(a,b).$ 

For any $y \in \pi^{(d-1)}(\nHat{S}_{\sigma}),$ we have that 
$$\pi_{d-1}(y, \nHat{S}_{\sigma}) = \{ (s_1=y, s_2, \dots, s_d) \in V_d \ |  \  s_k \in \Omega([0, \frac{z(\sigma,k)}{z(\sigma,k-1)}s_{k-1}]), \forall 2 \le k \le d\}$$
Thus, 
\begin{eqnarray*}
&& \vol_{\Lambda^1}(\pi_{d-1}(y, \nHat{S}_{\sigma})) \\
&=& \int_{s_2 \in \Omega([0, \frac{z(\sigma,2)}{z(\sigma,1)}y])} \left(\int_{s_3 \in \Omega([0, \frac{z(\sigma,3)}{z(\sigma,2)}s_{2}])} \! \! \! \! \! \cdots \left(\int_{s_d \in \Omega([0, \frac{z(\sigma,d)}{z(\sigma,d-1)}s_{d-1}])} \! \! \! \! \! 1 \ d s_d\right) \cdots d s_3\right) d s_2 \\
&=& |y|^{d-1} \prod_{j=2}^d \frac{1}{d+1-j} \left| \frac{z(\sigma,j)}{z(\sigma,j-1)} \right|^{d+1-j} 
= \frac{|y|^{d-1}}{(d-1)!} \left| \frac{\prod_{j=2}^d z(\sigma,j)}{z(\sigma,1)^{d-1}} \right|.
\end{eqnarray*}
It is clear that $\pi^{(d-1)}(\nHat{S}_{\sigma}) = \Omega([0, z(\sigma,1)]).$ Therefore, for any $y \in \pi^{(d-1)}(\nHat{S}_{\sigma})$, we have that $\sgn(y) = \sgn(z(\sigma,1)).$

Therefore, we have that 
$$\vol_{\Lambda^1}(\pi_{d-1}(y, \nHat{S}_{\sigma})) = \frac{y^{d-1}}{(d-1)!} \frac{\left| \prod_{j=2}^d z(\sigma,j) \right|}{z(\sigma,1)^{d-1}}.$$
Compare this with \eqref{svolnHat}, we see that it is left to show that
$$\sum_{y \in \Lambda_1 \cap \Omega([0, z(\sigma,1)]} y^{d-1} = \sgn(z(\sigma,1)) P_{d-1}(z(\sigma,1)).$$

However, since $P$ is an integral polytope, we have that $z(\sigma,1) = x_{d+1,1} - x_{\sigma(1),1}$ is an integer. If $z(\sigma,1) > 0,$ then 
$$\sum_{y \in \Lambda_1 \cap \Omega([0, z(\sigma,1)]} y^{d-1} = \sum_{0< y \le z(\sigma,1)} y^{d-1} = \sum_{y=1}^{z(\sigma,1)} y^{d-1} = P_{d-1}(z(\sigma,1)).$$
If $z(\sigma,1)<0,$ then
\begin{eqnarray*}
& &\sum_{y \in \Lambda_1 \cap \Omega([0, z(\sigma,1)]} y^{d-1} = \sum_{z(\sigma,1)< y \le 0} y^{d-1} = (-1)^{d-1} \sum_{0 \le y < -z(\sigma,1)} y^{d-1}\\
&=&  (-1)^{d-1} \sum_{y=1}^{-z(\sigma,1)-1} y^{d-1} = (-1)^{d-1} P_{d-1}(-z(\sigma,1)-1) = - P_{d-1}(z(\sigma,1)),
\end{eqnarray*}
where the last equality follows from \eqref{recp}.
\end{proof}

\subsection{Final calculation}
By \eqref{svol2Omega}, \eqref{fdecomp}, Lemma \ref{2hat}, and \eqref{svolnHat}, we have that
$$\svol_{\Gamma}^1(P) = \sum_{\sigma \in \S_d} \sgn(\sigma, P) \frac{1}{(d-1)!} \frac{\left|\prod_{j=1}^d z(\sigma,j)\right|}{z(\sigma,1)^{d}} P_{d-1}(z(\sigma,1)),$$
where $\sgn(\sigma, P)$ is defined as in \eqref{sign}.
Observing that
$$\sgn(\det(X(\sigma,d))) = \sgn(\sigma) \sgn(\det(X(\1,d))),$$ 
we conclude that to prove Proposition \ref{k=1simplex}, it suffices to prove the following lemma.

\begin{lem}\label{k=1simplex1}
Let $P \subset V = V_d$ be a $d$-dimensional ($0$-)integral simplex in fully general position with vertex set
$\vert = \{v_1, v_2, \dots, v_{d+1}\},$ where the coordinates of $v_i$ (with respect to $\e$)
are $ \x_i = (x_{i,1}, x_{i,2}, \dots, x_{i,d}).$ Then
\begin{equation}\label{signedsum}
\sum_{\sigma \in \S_d} \sgn(\sigma) \frac{1}{(d-1)!} \frac{\prod_{j=1}^d z(\sigma,j)}{z(\sigma,1)^{d}} P_{d-1}(z(\sigma,1)) =\frac{1}{d!}\det(X(\1, d)).
\end{equation}
\end{lem}

We need two more results from \cite{lattice-face}.
In Subsection 7.1 of \cite{lattice-face}, it was proved that 
\begin{equation}\label{rewritedet}
\frac{1}{d!}\det(X(\1, d)) = \frac{1}{d!} \sum_{\sigma \in \S_d} \sgn(\sigma) \prod_{j=1}^d z(\sigma, j).
\end{equation}
We also have the following lemma.
 
\begin{lem}[Lemma 7.3 in \cite{lattice-face}]\label{zerofive}
For any nonnegative integers $\ell,m$ satisfying $0 \le \ell +  m \le d-2,$  given $p(y_1, \dots, y_{\ell})$
a function on $\ell$ variables, let $q(\sigma) = p(z(\sigma,1),
\dots, z(\sigma,\ell)), \forall \sigma \in \S_d.$ Then
\begin{equation}\label{zero5}
\sum_{\sigma \in \S_d} \sgn(\sigma) q(\sigma) \frac{\prod_{j=\ell +
1}^{d} {z}(\sigma,j) }{({z}(\sigma,\ell+1))^{m+1}} = 0.
\end{equation}
\end{lem}

\begin{proof}[Proof of Lemma \ref{k=1simplex1}] 
By \eqref{rewritedet}, we only need to show the left hand side of \eqref{signedsum} is equal to $\frac{1}{d!} \sum_{\sigma \in \S_d} \sgn(\sigma) \prod_{j=1}^d z(\sigma, j)$ as well.

By \eqref{pdeg}, \eqref{pconst} and \eqref{plead}, the $(d-1)$th power sum polynomial can be expressed as $$P_{d-1}(x) = \frac{1}{d} x^d + \sum_{i=1}^{d-1} c_i x^i,$$
for some real numbers $c_i$'s. But for $i: 1 \le i \le d-1,$
$$\! \sum_{\sigma \in \S_d} \sgn(\sigma) \frac{1}{(d-1)!} \frac{\prod_{j=1}^d z(\sigma,j)}{z(\sigma,1)^{d}} c_i z(\sigma,1)^i=\! \frac{c_i}{(d-1)!} \sum_{\sigma \in \S_d} \sgn(\sigma) \frac{\prod_{j=1}^d z(\sigma,j)}{z(\sigma,1)^{d-i}} = 0,$$
where the second equality follows from Lemma \ref{zerofive} with $\ell = 0,$ $m = d-i-1,$ and $p =1.$ Hence, the left hand side of \eqref{signedsum} becomes
$$\sum_{\sigma \in \S_d} \sgn(\sigma) \frac{1}{(d-1)!} \frac{\prod_{j=1}^d z(\sigma,j)}{z(\sigma,1)^{d}} \frac{1}{d} z(\sigma,1)^d = \frac{1}{d!} \sum_{\sigma \in \S_d} \sgn(\sigma) \prod_{j=1}^d z(\sigma, j).$$
\end{proof}

We conclude by remarking that it is possible to modify our arguments in this section to work for general $k.$ In other words, with a modified argument, one can prove that Proposition \ref{mainVol1} is true for a simplex $P$ without restricting to $k=1$. However, since the modified argument is more complicated, and more importantly, unlike the case when $k=1,$ we cannot reduce the general case to the simplex case (which we will discuss later in Remark \ref{bigkFail}), we do not include them here.

\section{Properties of slices and projections}\label{propertySP}
Before we move on to finishing the proof of Proposition \ref{mainVol1} and thus Theorem \ref{mainVol}, we will discuss properties of slices and projections of certain polytopes. The results presented in this section will be used both in the proof of Proposition \ref{mainVol1} in Section \ref{pfVol} and in the proof of Theorem \ref{mainEhr} in Section \ref{pfEhr}.

Let $P$ be a polytope in $V$. One sees that the projection $\pi^{(D-m)}(P)$ is a polytope in $V_{m},$ and the slice $\pi_{D-m}(\y, P)$ is a polytope in $\pi_{D-m}(\y, V) \cong V^m,$ for any $0 \le m \le D$ and $\y \in \pi^{(D-m)}(P).$ We will discuss properties of projections and slices of $P$ under various conditions.

We denote by $\partial P$ and $\Int(P)$ the relative boundary and the relative interior of $P.$ Since in this paper we only discuss the relative boundary and the relative interior of polytopes, we will omit ``relative'' and just say ``boundary'' and ``interior''. We denote by $\F(P)$ and $\F_{\ge \ell}(P)$ the set of faces and the set of faces of dimension at least $\ell$, respectively, of $P.$  

Let $Y \subset V_{m},$ and $S \subset V.$ We define 
$$\pi_{D-m}(Y, S) = \bigcup_{\y \in Y} \pi_{D-m}(\y, S).$$ 

\begin{lem}\label{isom}
Suppose $0 < \ell \le m \le D$ and $P \subset V$ is a polytope whose $\ell$-faces are all in affinely general position. For any $(\ell-1)$-face $F_0$ of $\pi^{(D-m)}(P),$ we have that $F: = \pi_{D-m}(F_0, P)$ is an $(\ell-1)$-face of $P.$

Therefore, $\pi_{D-m}(\y, P)$ is one point on the boundary of $P$, for any point $\y$ of an $(\ell-1)$-face of $\pi^{(D-m)}(P).$ 
\end{lem}
\begin{proof}
Since $F_0$ is a face of  $\pi^{(D-m)}(P),$ there exists an $(m-1)$-dimensional affine space $H_0$ in $V_m$ such that $\pi^{(D-m)}(P)$ is on one side of $H_0$ and $F_0 = H_0 \cap \pi^{(D-m)}(P).$ Let $H$ be the inverse image of $H_0$ under $\pi^{(D-m)}.$ One sees that $H$ is a $(D-1)$-dimensional affine space in $V$, of which $P$ is on one side. It is clear that $F = \pi_{D-m}(F_0, P) = H \cap P,$ thus is a face of $P.$ Because $\pi^{(D-m)}(F) = F_0,$ we must have that $\dim(F) \ge \dim(F_0) = \ell-1.$ Suppose $\dim(F) \ge \ell.$ Let $F'$ be an $\ell$-face of $F.$ Then $F'$ is a face of $P,$ thus is in affinely general position. Therefore, by Corollary \ref{propIntGen}/(iii), $\dim(\pi^{(D-m)}(F')) = \ell.$ However, $\pi^{(D-m)}(F') \subset \pi^{(D-m)}(F) = F_0,$ which contradicts the hypothesis that $F_0$ has dimension $\ell-1.$
\end{proof}

\begin{cor}\label{charProj}
Suppose $0 \le k \le d \le D$ and $P \subset V$ is a $k$-integral polytope. Then $\pi^{(D-k)}(P)$ is a $k$-dimensional fully integral polytope. 
\end{cor}

\begin{proof}
By Lemma \ref{charRank}/(i), $\pi^{(D-k)}(P)$ is $k$-dimensional. The only $k$-face of $\pi^{(D-k)}(P)$ is itself. Since $\aff(\pi^{(D-k)}(P)) = V_k,$ it is automatically integral. For any $\ell < k,$  and any $\ell$-face $F_0$ of $\pi^{(D-k)}(P),$ by Lemma \ref{isom}, there exists an $\ell$-face $F$ of $P,$ such that $F = \pi_{D-k}(F_0, P).$ Thus, $\pi^{(D-k)}(F) = F_0.$ However, $F$ is affinely integral, which means $\aff(F)$ is integral. By Corollary \ref{propIntGen}/(ii), $\pi^{(D-k)}(\aff(F))$ is integral. One checks that $\aff(F_0) = \pi^{(D-k)}(\aff(F))$. Hence, $F_0$ is affinely integral. Therefore, $\pi^{(D-k)}(P)$ is fully integral.
\end{proof}

\begin{lem}\label{charDimSliceGen}
Suppose $0 \le k \le d \le D$ and $P \subset V$ is a $d$-dimensional polytope whose $k$-faces are all in affinely general position. Let $\y \in \pi^{(D-k)}(P).$ 
\begin{ilist}
\itm If $\y$ is an interior point of $\pi^{(D-k)}(P)$, then the slice $\pi_{D-k}(\y, P)$ is a $(d-k)$-dimensional polytope containing an interior point of $P.$
\itm If $\y$ is on the boundary of $\pi^{(D-k)}(P)$, then the slice $\pi_{D-k}(\y, P)$ is one point on the boundary of $P.$ (Note that this case would not happen if $k = 0$ because $\pi^{(D-k)}(P) = V_0$ does not have boundary points.)
\end{ilist}
\end{lem}

\begin{proof}
Let $U = \aff(P).$ Because $k$-faces of $P$ are in affinely general position, one checks that $\pi^{(D-k)}(U) = V_k.$ Therefore, $\pi_{D-k}(\y, U)$ is $(d-k)$-dimensional, for any $\y \in V_k.$  

Since $\pi^{(D-k)}$ is a continuous open map and $P$ is convex, we have that 
$$\pi^{(D-k)}(\Int(P)) = \Int(\pi^{(D-k)}(P)).$$
Thus,  if $\y$ is an interior point of $\pi^{(D-k)}(P),$ then $\pi_{D-k}(\y, P)$ contains an interior point of $P$ and $\dim(\pi_{D-k}(\y, P)) = \dim(\pi_{D-k}(\y, U)) = d-k.$ Now suppose $\y$ is a boundary point of $\pi^{(D-k)}(P).$ Then $\y$ is on a facet, which is a $(k-1)$-face, of $\pi^{(D-k)}(P).$ Hence, (ii) follows immediately from Lemma \ref{isom} with $\ell = m = k.$
\end{proof}

\begin{cor}\label{restrict2int}
Suppose $k < d$ and $P \subset V$ is a $d$-dimensional rational polytope.
\begin{ilist}
\itm  If all of the $k$-faces of $P$ are affinely integral, then
\begin{equation*}
\sum_{\y \in \pi^{(D-k)}(\Lambda_P)} \! \! \! \! \! \vol_{\Lambda^k \cap \Lambda_{\lin(P)}}(\pi_{D-k}(\y,P)) = \sum_{\y \in \Lambda \cap \Int(\pi^{(D-k)}(P))} \! \! \! \! \! \vol_{\Lambda^k \cap \Lambda_{\lin(P)}}(\pi_{D-k}(\y,P)).
\end{equation*}
\itm If $P$ is central and all of the $k$-faces of $P$ are in affinely general position, then
$$\svol_{\Lambda_P}^k(P) = \sum_{\y \in \pi^{(D-k)}(\Lambda_P)  \cap \Int(\pi^{(D-k)}(P))} \vol_{\Lambda^k \cap \Lambda_P}(\pi_{D-k}(\y,P)).$$ 

\end{ilist}
\end{cor}

\begin{proof}
We prove (ii) first. 
Since $d-k >0$, by Lemma \ref{charRank}/(ii) and Lemma \ref{charDimSliceGen}, for any $\y$ on the boundary of $\pi^{(D-k)}(P),$ we have that $\vol_{\Lambda^k \cap \Lambda_P}(\pi_{D-k}(\y,P)) = 0.$ Thus, (ii) follows.

Now we prove (i). One checks that it is enough to prove (i) under the assumption that $P$ is central. Then (i) follows from (ii), Corollary \ref{propIntGen}/(i) and Lemma \ref{charRank}/(i).
\end{proof}

The polytope $\pi_{D-k}(\y, P)$ is the intersection of the polytope $P$ with the affine space $\pi_{D-k}(\y, V).$ Hence, faces of $\pi_{D-k}(\y, P)$ are exactly the intersections of faces of $P$ with $\pi_{D-k}(\y, V).$ In other words,
\begin{equation}\label{relfaceofslice}
\F(\pi_{D-k}(\y, P)) = \{ \pi_{D-k}(\y, F) \ | \ F \in \F(P) \}.
\end{equation}
The following lemma discuss details of the relation between $\F(\pi_{D-k}(\y, P))$ and $\F(P).$

\begin{lem}\label{cnntFaceOfSlice}
Suppose $P \subset V$ is a $d$-dimensional polytope whose $k$-faces are all in affinely general position. Let $\y \in \pi^{(D-k)}(P).$ For any nonempty face $F_0$ of $\pi_{D-k}(\y,P),$ there exists a face $F$ of $P$ such that
\begin{ilist}
\itm $\dim(F) = \dim(F_0) + k,$
\itm $F_0 = \pi_{D-k}(\y, F).$
\end{ilist}

Furthermore, if $\dim(F_0) > 0,$ then such an $F$ is unique and we also have that $\y \in \Int(\pi^{(D-k)}(F)).$

\end{lem}

\begin{proof}
Suppose $\dim(F_0) = r.$ 
Because of \eqref{relfaceofslice}, there exists a face $F$ of $P$ such that $F_0 = \pi_{D-k}(\y, F).$ 
We consider 
\begin{equation}\label{psbFace}
\{ F \in \F(P) \ | \ F_0 = \pi_{D-k}(\y, F)\}.
\end{equation}

We claim that if \eqref{psbFace} contains one face of dimension less than or equal to $k,$ then $\dim(F_0) = 0$ and \eqref{psbFace} contains one face of dimension $k.$ Suppose $F$ is a face in \eqref{psbFace} of dimension less than or equal to $k,$ then $F$ is contained in a $k$-face $F'$ of $P.$ Because $F'$ is in affinely general position, it follows from Corollary \ref{propIntGen}/(iii) that $\pi_{D-k}(\y, F')$ contains at most one point. However, $F_0 = \pi_{D-k}(\y, F) \subset \pi_{D-k}(\y, F')$ and $F_0$ contains at least one point.
We must have that $F_0 = \pi_{D-k}(\y, F')$ is one point and $F'$ is in \eqref{psbFace}.

Suppose $\dim(F_0)=0,$ i.e., $F_0$ is a point. We want to show that the set \eqref{psbFace} contains one face of dimension $k.$ We have already seen that this is true if \eqref{psbFace} contains one face $F$ of dimension less than or equal to $k.$ Suppose all the faces in \eqref{psbFace} have dimension greater than $k,$ and let $F$ be one that has the smallest dimension. All the $k$-faces of $F$ are $k$-faces of $P,$ thus are in affinely general position. Hence, by Lemma \ref{charDimSliceGen}, we must have that $F_0 = \pi_{D-k}(\y, F)$ is on the boundary of $F.$ Therefore, there exists a proper face $F'$ of $F$ such that $F_0 \in F'.$ It is clear that $F_0 = \pi_{D-k}(\y, F').$ However, $\dim(F')<\dim(F),$ which contradicts the choice of $F.$ 

Suppose $\dim(F_0)>0.$ We want to show that the set \eqref{psbFace} contains exactly one face, which has dimension $\dim(F_0) + k.$ Let $F$ be a face in \eqref{psbFace}. By our claim, $\dim(F)>k.$ Again, all the $k$-faces of $F$ are in affinely general position. Thus, by Lemma \ref{charDimSliceGen}, we have that $\y \in \Int(\pi^{(D-k)}(F))$ and $F_0 = \pi_{D-k}(\y, F)$ has dimension $\dim(F)-k.$ Hence, $\dim(F) = \dim(F_0)+k.$ The uniqueness follows from the observation that if $F_1, F_2 \in \eqref{psbFace},$ then so is $F_1 \cap F_2.$
\end{proof}

\begin{cor}\label{bij}
Suppose $P \subset V$ is a $d$-dimensional polytope whose $k$-faces are all in affinely general position. Let $\y \in \Int(\pi^{(D-k)}(P)).$ 
\begin{ilist}
\itm Suppose $1 \le r \le d-k.$ The map $\pi_{D-k}(\y, \cdot)$ gives a bijection between the set of $(r+k)$ faces of $P$ satisfying $\y \in \Int(\pi^{(D-k)}(F))$ and the set of $r$-faces of $\pi_{D-k}(\y,P).$

\itm The map $\pi_{D-k}(\y, \cdot)$ gives a bijection between the set $\{ F \in \F_{\ge k+1}(P) \ | \ y \in \Int(\pi^{(D-k)}(F)) \}$ and the set $\F_{\ge 1}(\pi_{D-k}(\y,P)).$

\itm Suppose $P$ is a rational polytope. If $F$ is a face of $P$ that corresponds under $\pi_{D-k}(\y, \cdot)$ to a face $F_0$ of $\pi_{D-k}(\y,P),$ then $\Lambda_{\lin(F_0)} = \Lambda^k \cap \Lambda_{\lin(F)}.$
\end{ilist}
\end{cor}

\begin{proof}
(i) follows from Lemma \ref{cnntFaceOfSlice} and Lemma \ref{charDimSliceGen}/(i). (ii) follows from (i) and the fact that $P$ is $d$-dimensional. 

If $F$ is in bijection with $F_0$ under $\pi_{D-k}(\y, \cdot),$ then $F_0 = \pi_{D-k}(\y, F).$ One checks that $\lin(F_0) = V^k \cap \lin(F).$ Then (iii) follows.
\end{proof}

The original definitions of general position and integrality of an affine space require the space to sit inside the space spanned by the given lattice. It cannot be applied to slices of a polytope $P$ since they are in the affine spaces $\pi_{D-k}(\y, V).$ Hence, we extend these definitions. 


\begin{defn}
Let $U$ be an $\ell$-dimensional rational subspace of $V.$ Suppose $\f = (\bf^1, \dots, \bf^\ell)$ is a basis of the lattice $U \cap \Lambda.$ 
Let $W$ be an affine space of $V$ such that $\lin(W)$ is a subspace of $U.$ 
\begin{alist}
\itm 
We say $W$ is {\it integral with respect to $\f$} if $W$ contains a lattice point and $\lin(W)$ is integral with respect to $\f.$ 
\itm 
We say $W$ is {\it in general position with respect to $\f$} if  $\lin(W)$ is in general position with respect to $\f.$
\end{alist}
Similarly, we extend the definitions of {\it affinely integral, $k$-integral, affinely general position, and $k$-general position with respect to $\f$} for polytopes.
\end{defn}

\begin{lem}\label{sliceGenInt1}
Let $W$ be an $\ell$-dimensional affine space of $V.$ Suppose $m: 0 \le m \le \ell.$
\begin{ilist}
\itm If $W$ is integral, then $\pi_{D-m}(\y, W)$ is 
integral with respect to $(\be^{m+1}, \dots, \be^{D})$ for any $\y \in \Lambda_m.$
\itm If $W$ is in general position, then $\pi_{D-m}(\y, W)$ is 
in general position with respect to $(\be^{m+1}, \dots, \be^{D})$ for any $\y \in V_m.$
\end{ilist}
\end{lem}

\begin{proof} 
We will only prove (i), as (ii) can be proved similarly.
By Lemma \ref{charInt}, $\lin(W)$ is integral and $W = \lin(W) + \z$ for some $\z \in \Lambda \cap W.$ Then $\lin(W)$ is the row space of a matrix of the form $(I, J)$ where $I$ is the $\ell \times \ell$ identity matrix and $J$ is an $\ell \times (D-\ell)$ integer matrix.  Let $A$ be the top $m \times D$ submatrix of $(I,J)$ and $B$ be the bottom $(\ell-m) \times D$ submatrix of $(I,J).$ It is easy to verify that 
$$\pi_{D-m}(\y, W) =  \mbox{the row space of $B$} + ((\y- \pi^{(D-m)}(\z)) A + \z).$$
One sees that 
$$\lin(\pi_{D-m}(\y, W)) =  \mbox{the row space of $B$} \subset V^m.$$

Note that the first $m$ columns of $B$ are zeros, which is followed by the $(\ell-m)\times(\ell-m)$ identity matrix and an $(\ell-m)\times (D-\ell)$ integer matrix. Thus, $\lin(\pi_{D-m}(\y, W))$ is integral with respect to $(\be^{m+1}, \dots, \be^{D}).$ Since $(\y- \pi^{(D-m)}(\z)) A + \z$ is a lattice point in $W$, we conclude that $\pi_{D-m}(\y, W)$ is integral with respect to $(\be^{m+1}, \dots, \be^{D}).$ 
\end{proof}


\begin{lem}\label{sliceGenInt2}
Suppose $0 \le m \le \ell \le d$ and $P \subset V$ is a $d$-dimensional polytope. 
\begin{ilist}
\itm Suppose all the $\ell$-faces of $P$ are affinely integral. Then for each $\y$ in $\pi^{(D-m)}(P) \cap \Lambda,$ we have that every $(\ell-m)$-face of $\pi_{D-m}(\y, P)$ is affinely integral with respect to $(\be^{m+1}, \dots, \be^D)$.
\itm Suppose all the $\ell$-faces of $P$ are in affinely general position. Then for each $\y$ in $\pi^{(D-m)}(P)$ we have that every $(\ell-m)$-face of $\pi_{D-m}(\y, P)$ is in affinely general position with respect to $(\be^{m+1}, \dots, \be^D)$.
\end{ilist}
\end{lem}

\begin{proof} We only prove (i), and (ii) can be proved similarly.
Let $F_0$ be a face of $\pi_{D-m}(P)$ of dimension $\ell-m.$ By Lemma \ref{cnntFaceOfSlice}, there exists a face $F$ of $P$ such that $\dim(F) = \dim(F_0) + m = \ell$ and $F_0 = \pi_{D-m}(\y,  F).$ We have that $F$ is affinely integral, i.e., $\aff(F)$ is integral.  One checks that $\aff(F_0) = \aff(\pi_{D-m}(\y, F)) =\pi_{D-m}(\y, \aff(F)),$ which by Lemma \ref{sliceGenInt1}/(i) is integral with respect to $(\be^{m+1}, \dots, \be^D)$.
Therefore, $F_0$ is affinely integral with respect to $(\be^{m+1}, \dots, \be^D)$.
\end{proof}



\begin{cor}\label{sliceInt}
 Suppose all the $k$-faces of $P$ are affinely integral. Then for each $\y \in \pi^{(D-k)}(P) \cap \Lambda,$ we have that $\pi_{D-k}(\y, P)$ is an integral polytope.
\end{cor}
\begin{proof}
This follows from Lemma \ref{sliceGenInt2}/(i) with $m = \ell = k.$
\end{proof}

\section{Proof of Proposition \ref{mainVol1} and Theorem \ref{mainVol}}\label{pfVol}
We will first show that Proposition \ref{mainVol1} is true when $k=1,$ and then prove it for any $k$ by induction.

\begin{prop}\label{k=1}
Proposition \ref{mainVol1} is true when $k=1.$
\end{prop}

\begin{lem}\label{trisvol}
Suppose $0<k<d.$ Let $P$ be a $d$-dimensional central rational polytope in $V.$ Suppose $P$ has a triangulation $\bigsqcup Q_i$ consisting of  simplices in $k$-general position. Let $\Gamma = \Lambda_P.$ It is clear that $\aff(P) = \aff(Q_i).$ So $Q_i$ is central and $\Gamma = \Lambda_{Q_i}.$ Then
$$\svol_{\Gamma}^k(P) = \sum \svol_{\Gamma}^k(Q_i).$$
\end{lem}

\begin{proof}
Because $\bigsqcup Q_i$ is a triangulation of $P,$ we have that
$$[P] = \sum_{i} [Q_i] + \sum_{j} \epsilon_i [F_j],$$
where $\{ F_j \}$ is a collection of proper faces of $Q_i$'s, and $\epsilon_i = \pm 1.$ It is enough to show that 
\begin{equation}\label{facevanish}
\svol_{\Gamma}^k(F_j) = 0, \forall j.
\end{equation}
By Lemma \ref{charRank}/(ii), we have that $\Lambda^k \cap \Gamma$ is of rank $d-k.$ Since $F_j$ is a proper face of some $Q_i,$ we have that $F_j$ is contained in a facet $F_j'$ of $P.$ The facet $F_j'$ has dimension $d-1$ and is in $k$-general position. It follows from Lemma \ref{charDimSliceGen} that for each $\y \in \pi^{(D-k)}(F_j'),$ the dimension of $\pi_{D-k}(\y, F_j')$ is less than $d-k.$ Since $\pi_{D-k}(\y, F_j)$ is a subset of $\pi_{D-k}(\y, F_j'),$ its dimension is less than $d-k$ as well. Hence, $\vol_{\Lambda^k \cap \Gamma}(\pi_{D-k}(\y, F_j)) = 0.$ \eqref{facevanish} follows.
\end{proof}

\begin{lem}\label{triGen}
Let $P$ be a $d$-dimensional polytope in $1$-general position. Then there exists a triangulation $\bigsqcup Q_i$ of $P$ without introducing new vertices  such that each $Q_i$ is a simplex in $1$-general position.
\end{lem}

\begin{proof}
We prove the lemma by induction on $d.$ If $d=1,$ $P$ itself is a simplex, thus the lemma holds. 
Suppose now the lemma is true for polytopes of dimension less than $d.$
Let $y_0$ be the smallest number in the closed interval $\pi^{(D-1)}(P).$ By Lemma \ref{isom} with $\ell = m =1$, we have that $v_0 := \pi_{D-1}(y_0, P)$ is a vertex of $P.$ Clearly, $v_0$ is the only vertex of $P$ whose first coordinate is $y_0.$ Let $\{ F_i\}$ be the set of facets of $P$ that do not contain $v_0.$ Each $F_i$ is in $1$-general position. Thus, by the induction hypothesis, there exists a triangulation $\bigsqcup_{j=1}^{\ell_i} Q_{i,j}$ of $F_i$ without introducing new vertices  such that each $Q_{i,j}$ is a $(d-1)$-dimensional simplex in $1$-general position.
One checks that $$\bigsqcup_{i} \bigsqcup_{j=1}^{\ell_i} \conv(v_0, Q_{i,j})$$ is a triangulation of $P$ without new vertices. For any edge $e$ in a $d$-dimensional simplex $\conv(v_0, Q_{i,j})$ involved in the triangulation, $e$ is either an edge in $Q_{i,j}$ or is an edge connecting $v_0$ and a vertex $v$ in $Q_{i,j}.$ For the former case, since $Q_{i,j}$ is in $1$-general position, we have that $\aff(e)$ is in general position. For the latter case, because $\bigsqcup_{j=1}^{\ell_i} Q_{i,j}$ does not introduce new vertices, $v$ is a vertex of $F_i,$ thus is a vertex of $P.$ Hence, $v$ has different first coordinate from $v_0$ and $\aff(e)$ is also in general position. Therefore, $\conv(v_0, Q_{i,j})$ is in $1$-general position.
\end{proof}

\begin{proof}[Proof of Proposition \ref{k=1}]
By Lemma \ref{triGen}, there exists a triangulation $\bigsqcup Q_i$ of $P$ without introducing new vertices such that each $Q_i$ is a simplex in $1$-general position. $P$ is ($0$-)integral, so all the vertices of $P$ are integral. Since the triangulation does not introduce new vertices, each $Q_i$ is ($0$-)integral. Therefore, the proposition follows from Proposition \ref{k=1simplex}, Lemma \ref{trisvol} and the fact that $\vol_{\Gamma}(P) = \sum \vol_{\Gamma}(Q_i),$ where $\Gamma = \Lambda_P = \Lambda_{Q_i}.$
\end{proof}

\begin{rem}\label{bigkFail}
Lemma \ref{triGen} is essential for proving Proposition \ref{mainVol1} for $k=1.$ It reduces the problem to the case when $P$ is a simplex. However, we are unable to find a lemma that works for larger $k.$ Indeed, if we want to reduce the problem of Proposition \ref{mainVol1} for $k >1$ to the simplex case, we need to show that there exists a triangulation of $P$ such that each simplex is $(k-1)$-integral and in $k$-general position. While being in $k$-general position is easy to achieve, it is not clear how we can guarantee $(k-1)$-integral.

Therefore, for $k>1,$ we will prove Proposition \ref{mainVol1} differently, using induction on $k.$
\end{rem}

\begin{proof}[Proof of Proposition \ref{mainVol1} (and Theorem \ref{mainVol})]
We will prove the proposition by induction on $k.$ The base case $k=1$ follows from Proposition \ref{k=1}. Let $k_0: 2 \le k_0 <d.$ Now we assume that the theorem holds for any $k \le k_0,$ and we will prove the theorem for $k =k_0.$ By Corollary \ref{restrict2int}/(ii), it is sufficient to show that
\begin{equation}\label{sumIntFor}
\vol_{\Lambda_P}(P) = \sum_{\z \in \pi^{(D-k)}(\Lambda_P)   \cap \Int(\pi^{(D-k_0)}(P))} \vol_{\Lambda^{k_0} \cap \Lambda_P}(\pi_{D-k_0}(\z,P)).
\end{equation}
Since $k_0 \ge 2,$ we have that $P$ is an ($0$-)integral polytope in $1$-general position. Hence, by Proposition \ref{k=1} and Corollary \ref{restrict2int}/(ii), we have that 
\begin{eqnarray*}
\vol_{\Lambda_P}(P) &=& \sum_{y \in \pi^{(D-1)}(\Lambda_P) \cap \Int(\pi^{(D-1)}(P))} \vol_{\Lambda^1 \cap \Lambda_P}(\pi_{D-1}(y,P))
\end{eqnarray*}
It follows from Lemma \ref{sliceGenInt2} that the slice $\pi_{D-1}(y,P)$ is $(k_0-2)$-integral and in $(k_0-1)$-general position with respect to $(\be^2, \dots, \be^D),$ for any $y \in \pi^{(D-k)}(\Lambda_P)  \cap \Int(\pi^{(D-1)}(P)).$ Thus, $\pi_{D-1}(y,P) - y$ is a central polytope in $V^1$ that is $(k_0-2)$-integral and in $(k_0-1)$-general position with respect to $(\be^2, \dots, \be^D).$ Moreover, the lattice in $V^1$ is $\Lambda^1,$ and the lattice we use to calculate the normalized volume of $\pi_{D-1}(y,P)- y$ is 
$$\Lambda^1 \cap \lin(\pi_{D-1}(y,P)- y) = \Lambda^1 \cap \pi_{D-1}(0, \lin(P)) = \Lambda^1 \cap (V^1 \cap \lin(P)) = \Lambda^1 \cap \Lambda_P.$$
Hence, we can apply the induction hypothesis to $\pi_{D-1}(y,P) - y$ and translate the result to $\pi_{D-1}(y,P).$ We then get
\begin{equation*}
\vol_{\Lambda^1 \cap \Lambda_P}(\pi_{D-1}(y,P)) = \sum_{\substack{\z \in \pi^{(D-k_0)}(\Lambda^1 \cap \Lambda_P) + y \\  \z \in \Int(\pi^{(D-k_0)}(\pi_{D-1}(y,P)))}} \! \! \! \vol_{\Lambda^{k_0} \cap \Lambda_P}(\pi_{D-k_0}(\z, \pi_{D-1}(y,P))).
\end{equation*}

One checks that for any $\z \in V_{k_0}$, we have that $\z \in \pi^{(D-k_0)}(\Lambda^1 \cap \Lambda_P) + y$ if and only if $\z \in \pi^{(D-k_0)}(\Lambda_P)$ and the first coordinate of $\z$ is $y.$ Furthermore, $\z \in  \Int(\pi^{(D-k_0)}(\pi_{D-1}(y,P))) = \Int(\pi_{D-1}(y, \pi^{(D-k_0)}(P)))$ if and only if the first coordinate of $\z$ is $y,$ and $\z \in \Int(\pi^{(D-k_0)}(P)).$ Also, if $\z \in  \pi^{(D-k_0)}(\pi_{D-1}(y,P))$, then $\pi_{D-k_0}(\z, \pi_{D-1}(y,P)) = \pi_{D-k_0}(\z, P).$
Therefore, we have
\begin{eqnarray*}
& &\vol_{\Lambda_P}(P) \\
&=&\! \! \! \! \! \! \! \! \! \! \sum_{y \in \pi^{(D-1)}(\Lambda_P) \cap \Int(\pi^{(D-1)}(P))}   \sum_{\substack{\z \in \pi^{(D-k_0)}(\Lambda^1 \cap \Lambda_P) + y \\  \z \in \Int(\pi^{(D-k_0)}(\pi_{D-1}(y,P)))}} \! \! \! \! \! \! \! \! \! \! \vol_{\Lambda^{k_0} \cap \Lambda_P}(\pi_{D-k_0}(\z, \pi_{D-1}(y,P)).\\
&=& \! \! \! \! \! \sum_{y \in \pi^{(D-1)}(\Lambda_P) \cap \Int(\pi^{(D-1)}(P))} \sum_{\substack{\z \in \pi^{(D-k_0)}(\Lambda_P) \cap \Int(\pi^{(D-k_0)}(P)) \\ \pi^{(D-1)}(\z) = y
}} \! \! \! \! \! \vol_{\Lambda^{k_0} \cap \Lambda_P}(\pi_{D-k_0}(\z, P)).
\end{eqnarray*}
Furthermore,
\begin{eqnarray*}
& &\pi^{(k_0-1)}(\Int(\pi^{(D-k_0)}(P))) \\
&=& \pi^{(k_0-1)}(\pi^{(D-k_0)}(\Int(P))) = \pi^{(D-1)}(\Int(P))=\Int(\pi^{(D-1)}(P)),
\end{eqnarray*}
which means that $\Int(\pi^{(D-1)}(P))$ contains all possible first coordinates of points in $\Int(\pi^{(D-k_0)}(P)).$ 
Therefore, \eqref{sumIntFor} follows.
\end{proof}

\section{Proof of Theorem \ref{mainEhr}}\label{pfEhr}

We first recall a result from \cite{berline-vergne}, which we will use to prove Theorem \ref{mainEhr}.

Let $P \subset V$ be a nonempty polyhedron and let $\v \in P$ be a point. We define the {\it tangent cone of $P$ at $\v$} by 
$$\tcone(P, \v) = \{ \v + \y \ | \ \v + \epsilon \y \in P \mbox{ for some $\epsilon > 0$} \}.$$ 
We define the {\it feasible cone of $P$ at $\v$} by
$$\fcone(P, \v) = \{ \y \ | \ \v + \epsilon \y \in P \mbox{ for some $\epsilon > 0$} \}.$$ 
Thus, $\tcone(P, \v) = \fcone(P, \v) + \v.$

For any face $F$ of $P$ and any two interior points $\v, \w$ of $F$, one checks that $\tcone(P, \v) = \tcone(P, \w)$ and $\fcone(P, \v) = \fcone(P, \w).$ Therefore, we can define the {\it tangent cone of $P$ along $F$} by 
$$\tcone(P, F) := \tcone(P, \v), \mbox{ for some (equivalently all) $\v \in \Int(F)$,}$$
and the {\it feasible cone of $P$ along $F$} by 
$$\fcone(P, F) := \fcone(P, \v), \mbox{ for some (equivalently all) $\v \in \Int(F)$.}$$

In \cite{berline-vergne}, Berline and Vergne showed that for any rational polytope $P,$
\begin{equation}\label{LPfromVol}
i(P) = \sum_{F \in \F_{\ge 0}(P)} \alpha(P, F) \vol_{\Lambda_{\lin(F)}}(F),
\end{equation}
where $\alpha(P,F)$ is a rational number which only depends on $\tcone(P, F).$ 
In this paper, we only need the results on $\alpha(P,F)$ when $P$ is an integral polytope, which we summarize in the following lemma:

\begin{lem}[Corollaries 30 and 17 in \cite{berline-vergne}]\label{int}
Suppose $P$ is an integral polytope. 
\begin{ilist}
\itm $\alpha(P,F)$ is determined by $\fcone(P,F).$ 
\itm $\sum_{\v \in \vert(P)} \alpha(P, \v) = 1.$
\end{ilist}
\end{lem}

We have the following corollary which immediately follows from \eqref{LPfromVol} and Lemma \ref{int}/(ii).
\begin{cor}\label{minus1}
If $P$ is an integral polytope, then
$$i(P)-1 = \sum_{F \in \F_{\ge 1}(P)} \alpha(P, F) \vol_{\Lambda_{\lin(F)}}(F).$$
Hence,
$$i(P, m)-1 = \sum_{F \in \F_{\ge 1}(P)} \alpha(P, F) \vol_{\Lambda_{\lin(F)}}(F) m^{\dim(F)}.$$
\end{cor}

\begin{lem}
Suppose $P \subset V$ is a $d$-dimensional polytope whose $k$-faces are all affinely integral. Let $F$ be a face of $P$ of dimension greater than $k.$ 
\begin{ilist}
\itm $\pi_{D-k}(\y, F)$ is a $(\dim(F)-k)$-face of $\pi_{D-k}(\y, P),$ for any lattice point $\y$ in $\Int(\pi^{(D-k)}(F)).$ 
\itm $\alpha(\pi_{D-k}(\y, P), \pi_{D-k}(\y, F))$ is invariant under different choices of lattice points $\y$ in  $\Int(\pi^{(D-k)}(P)).$
\end{ilist}
\end{lem}
\begin{proof}
(i) follows from Corollary \ref{bij}.
By Corollary \ref{sliceInt}, $\pi_{D-k}(\y, P)$ is integral. One checks that $\fcone(\pi_{D-k}(\y, P), \pi_{D-k}(\y, F))$ is exactly $\fcone(P,F) \cap V^k,$ which is independent of the choice of $\y.$ (ii) then follows from Lemma \ref{int}/(i).
\end{proof}

Given this lemma, we are able to give the following definition.

\begin{defn}
Suppose $P \subset V$ is a $d$-dimensional polytope whose $k$-faces are all affinely integral. Let $F$ be a face of $P$ of dimension greater than $k.$ We define
$$\beta_k(P, F) := \alpha(\pi_{D-k}(\y, P), \pi_{D-k}(\y, F)),$$ 
for some (equivalently all) $\y \in \Lambda \cap \Int(\pi^{(D-k)}(F))$.
\end{defn}






In the following lemma, we give formulas connecting numbers of lattice points of slices to $\beta_k(P,F)$'s and normalized volumes of faces $F$ of dimension at least $k+1.$

\begin{lem}\label{cnntMinus1}
Suppose $0 \le k < d \le D$ and $P  \subset V$ is a $k$-integral $d$-dimensional polytope. For any positive integer $m,$ we have that
\begin{equation}\label{sliceEhrFor}
\sum_{\y \in \Lambda \cap \pi^{(D-k)}(P)} (i(\pi_{D-k}(\y, P), m) - 1) = m^{-k} \! \! \! \! \! \sum_{F \in \F_{\ge k+1}(P)} \! \! \! \! \! \beta_k(P, F) \vol_{\Lambda_{\lin(F)}}(F) \ m^{\dim(F)}.
\end{equation}

Hence,  
\begin{equation}\label{chg2minus1}
\sum_{\y \in \Lambda \cap \pi^{(D-k)}(P)} (i(\pi_{D-k}(\y, P)) - 1) = \sum_{F \in \F_{\ge k+1}(P)} \beta_k(P, F) \vol_{\Lambda_{\lin(F)}}(F).
\end{equation}

In particular, if $k = d-1,$ then
\begin{equation}\label{d-1integral0}
\sum_{\y \in \Lambda \cap \pi^{(D-(d-1))}(P)} (i(\pi_{D-(d-1)}(\y, P)) - 1) = \vol_{\Lambda_{\lin(P)}}(P)
\end{equation}
\end{lem}

\begin{proof}
If $k = 0,$ then \eqref{sliceEhrFor} follows from Corollary \ref{minus1} and the observation that the left hand side of \eqref{sliceEhrFor} is $i(P,m)-1,$ and $\beta_0(P, F) = \alpha(P,F).$ Hence, we can assume $k \ge 1.$
By Corollary \ref{sliceInt}, $\pi_{D-k}(\y, P)$ is integral, for any $\y \in \Lambda \cap \pi^{(D-k)}(P).$ If $\y$ is on the boundary of $\pi^{(D-k)}(P),$ by Lemma \ref{charDimSliceGen}/(ii), $\pi_{D-k}(\y, P)$ is just one point. Therefore, $i(\pi_{D-k}(\y, P),m) - 1=0.$ Using this and Corollary \ref{minus1}, we have that the left hand side of \eqref{sliceEhrFor} is equal to
\begin{equation}\label{wrtSlice}
 \sum_{\y \in \Lambda \cap \Int(\pi^{(D-k)}(P))}  \sum_{F_0 \in \F_{\ge 1}(\pi_{D-k}(\y, P))} \alpha(\pi_{D-k}(\y, P), F_0) \vol_{\Lambda_{\lin(F_0)}}(F_0) \ m^{\dim(F_0)}.
\end{equation}
However, by Corollary \ref{bij}, we can rewrite \eqref{wrtSlice} as
\begin{eqnarray*}
& &\sum_{\y \in \Lambda \cap \Int(\pi^{(D-k)}(P))}  \sum_{\substack{F \in \F_{\ge k+1}(P)\\ \y \in \Int(\pi^{(D-k)}(F))}} \Bigl[\alpha(\pi_{D-k}(\y, P), \pi_{D-k}(\y, F)) \times \\
& & \hspace{5cm} \vol_{\Lambda^k \cap \Lambda_{\lin(F)}}(\pi_{D-k}(\y, F)) \  m^{\dim(F) - k} \Bigr] \\
&=& \sum_{\y \in \Lambda \cap \Int(\pi^{(D-k)}(P))}  \! \! \! \sum_{\substack{F \in \F_{\ge k+1}(P)\\ \y \in \Int(\pi^{(D-k)}(F))}} \! \! \! \! \! \! \beta_k(P, F) \vol_{\Lambda^k \cap \Lambda_{\lin(F)}}(\pi_{D-k}(\y, F)) \ m^{\dim(F)-k}
\end{eqnarray*}
However, for any face $F \in \F_{\ge k+1}(P),$ it follows from the hypothesis that $P$ is $k$-integral that  $F$ is $k$-integral as well. By Corollary \ref{charProj}, we have that $\dim(\pi^{(D-k)}(F)) = k = \dim(\pi^{(D-k)}(P)).$ Hence, $\Int(\pi^{(D-k)}(F)) \subset \Int(\pi^{(D-k)}(P)).$ Therefore, we can simplify the above formula further and conclude that
\begin{eqnarray*}
& &\sum_{\y \in \Lambda \cap \pi^{(D-k)}(P)} (i(\pi_{D-k}(\y, P), m) - 1) \\
&=& \sum_{\substack{F \in \F_{\ge k+1}(P)\\ \y \in \Lambda \cap \Int(\pi^{(D-k)}(F))}} \beta_k(P, F) \vol_{\Lambda^k \cap \Lambda_{\lin(F)}}(\pi_{D-k}(\y, F)) \ m^{\dim(F)-k} \\
&=& m^{-k} \sum_{F \in \F_{\ge k+1}(P)} \beta_k(P, F) \ m^{\dim(F)} \! \! \! \! \! \sum_{\y \in \Lambda \cap \Int(\pi^{(D-k)}(F))} \! \! \! \! \! \vol_{\Lambda^k \cap \Lambda_{\lin(F)}}(\pi_{D-k}(\y, F)).
\end{eqnarray*}
Since any face $F \in \F_{\ge k+1}(P)$ is $k$-integral, and any $k$-integral polytope is $(k-1)$-integral and in $k$-general position. Therefore, by Theorem \ref{mainVol} and Corollary \ref{restrict2int}/(i),  the identity \eqref{sliceEhrFor} holds.

\eqref{chg2minus1} is obtained from \eqref{sliceEhrFor} by plugging in $m=1.$

Finally, in the case of $k=d-1,$ note that the only face $F$ in $\F_{\ge d}(P)$ is $P$, so \eqref{d-1integral0} would follow from \eqref{chg2minus1} given the identity $$\beta_{d-1}(P,P) =\alpha(\pi_{D-(d-1)}(\y, P), \pi_{D-(d-1)}(\y, P)) = 1.$$
However, in the construction of $\alpha(P,F)$ in \cite{berline-vergne}, $\alpha(Q,Q)$ is always $1$ for any rational polytope. So we have the above identity.
\end{proof}

\begin{rem}
It is possible to give a proof of \eqref{d-1integral0} that makes no use of the results of \cite{berline-vergne}. In fact, we only need the simple fact that the normalized volume of any $1$-dimensional integral polytope is equal to the number of lattice points minus $1,$ and then we can prove \eqref{d-1integral0} with arguments similar to those we used in proving \eqref{sliceEhrFor} (but without involving $\beta_k(P,F)$).


\end{rem}

\begin{prop}\label{cnntLPwProj}
Suppose $0 \le k < d \le D$ and $P  \subset V$ is a $k$-integral $d$-dimensional polytope. Then 
\begin{equation}\label{LPfromProjVol}
i(P) = i(\pi^{(D-k)}(P)) + \sum_{F \in \F_{\ge k+1}(P)} \beta_k(P, F) \vol_{\Lambda_{\lin(F)}}(F).
\end{equation}

In particular, if $k = d-1,$ then
\begin{equation}\label{d-1integral}
i(P) = i(\pi^{(D-(d-1))}(P)) + \vol_{\Lambda_{\lin(P)}}(P).
\end{equation}
\end{prop}

\begin{proof}
It is clear that 
\begin{eqnarray}
i(P) &=& \sum_{\y \in \Lambda \cap \pi^{(D-k)}(P)} i(\pi_{D-k}(\y, P)) \nonumber \\
&=&  i(\pi^{(D-k)}(P)) +  \sum_{\y \in \Lambda \cap \pi^{(D-k)}(P)} (i(\pi_{D-k}(\y, P)) - 1).  \nonumber
\end{eqnarray}
Then the proposition follows from Lemma \ref{cnntMinus1}.
\end{proof}

Using \eqref{d-1integral}, we are able to generalize the results in \cite{cyclic, lattice-face, note-lattice-face} to fully integral polytopes:

\begin{cor}\label{fullyint}
Suppose $0 \le d \le D$ and $P  \subset V$ is a $d$-dimensional fully integral polytope. Then
\begin{equation}\label{LP=sumvol}
i(P) = \sum_{j=0}^d \vol_{\Lambda_j} (\pi^{(D-j)}(P)).
\end{equation}
Hence, the Ehrhart polynomial of $P$ is given by
\begin{equation}\label{coeff=vol}
i(P,m) = \sum_{j=0}^d \vol_{\Lambda_j} (\pi^{(D-j)}(P)) m^j.
\end{equation}
\end{cor}

\begin{proof}
It is clear that the corollary holds for $d = 0.$ Now we assume $d \ge 1.$
$P$ is fully integral, so $P$ is $(d-1)$-integral and $\aff(P)$ is integral. Therefore, by Proposition \ref{cnntLPwProj} and Corollary \ref{propIntGen}/(ii), we have that $i(P) = i(\pi^{(D-(d-1))}(P)) + \vol_{\Lambda_{\lin(P)}}(P),$ and $\pi^{(D-d)}$ induces a bijection between $\aff(P) \cap \Lambda$ and $\Lambda_d.$ Because $\aff(P)$ contains a lattice point, $\aff(P) \cap \Lambda$ and $\Lambda_{\lin(P)}$ only differ by a lattice point. Hence, $\pi^{(D-d)}$ induces a bijection between $\Lambda_{\lin(P)}$ and $\Lambda_d.$ Therefore, we conclude that
\begin{equation*}
i(P) = i(\pi^{(D-(d-1))}(P)) + \vol_{\Lambda_d}(\pi^{(D-d)}(P)).
\end{equation*}
However, it follows from Corollary \ref{charProj} that $\pi^{(D-(d-1))}(P)$ is a $(d-1)$-dimensional fully integral polytope. Therefore, we obtain \eqref{LP=sumvol} by recursively applying the above identity.

Finally, by Lemma \ref{preDilateInt}, the dilation $m P$ of $P$ is fully integral, for any positive integer $m.$ Thus, \eqref{coeff=vol} follows.
\end{proof}




\begin{proof}[Proof of Theorem \ref{mainEhr}]
If $k=d,$ the theorem follows from \eqref{coeff=vol}. Now we assume $k<d.$ Let $m$ be a positive integer. By Lemma \ref{preDilateInt}, $mP$ is a $k$-integral polytope. For any $F \in \F_{\ge k+1}(P),$ one checks that $\beta_k(mP, mF) = \beta_k(P,F)$ and $\lin(mF) = \lin(F).$ Hence, by Proposition \ref{cnntLPwProj}, we have that
\begin{eqnarray*}
i(P, m) &=& i(\pi^{(D-k)}(P), m) + \sum_{F \in \F_{\ge k+1}(P)} \beta_k(P, F) \vol_{\Lambda_{\lin(F)}}(mF) \\
&=& i(\pi^{(D-k)}(P), m) + \sum_{F \in \F_{\ge k+1}(P)} \beta_k(P, F) \vol_{\Lambda_{\lin(F)}}(F) \ m^{\dim(F)}.
\end{eqnarray*}

By Corollary \ref{charProj}, $\pi^{(D-k)}(P)$ is a $k$-dimensional fully integral polytope. Therefore, we conclude our theorem using \eqref{sliceEhrFor} and \eqref{coeff=vol}.
\end{proof}

\bibliographystyle{amsplain}
\bibliography{gen}

\end{document}